 \documentclass[11pt]{article}

 \usepackage{amsfonts,amsmath,amssymb,mathrsfs}
 \usepackage{color}
 \usepackage[svgnames]{xcolor}
 \usepackage{harvard}
 \usepackage{bm}

\allowdisplaybreaks[4]

 \newtheorem{theorem}{Theorem}[section]
 \newtheorem{lemma}[theorem]{Lemma}
 \newtheorem{corollary}[theorem]{Corollary}
 \newtheorem{proposition}[theorem]{Proposition}
 \newtheorem{remark}[theorem]{Remark}
 \newtheorem{example}[theorem]{Example}
 \newtheorem{condition}[theorem]{Condition}

 \def\blemma{\begin{lemma}\sl{}
 \def\elemma{\end{lemma}}}
 
 \def\bproposition{\begin{proposition}\sl{}
 \def\eproposition{\end{proposition}}}

 \def\benumerate{\begin{enumerate}}\def\eenumerate{\end{enumerate}}
 \def\bitemize{\begin{itemize}}\def\eitemize{\end{itemize}}

 \def\beqlb{\begin{eqnarray}}
 \def\eeqlb{\end{eqnarray}}
 \def\beqnn{\begin{eqnarray*}}
 \def\eeqnn{\end{eqnarray*}}

 \def\qed{\hfill$\Box$\medskip}
 \def\bproof{\begin{proof}}\def\eproof{\qed\end{proof}}

 \def\<{\langle}\def\>{\rangle}

 \def\mbb{\mathbb}
 \def\mbf{\mathbf}\def\mrm{\mathrm}

 \def\ar{\!\!\!\!&}

 \def\d{\mrm{d}}\def\e{\mrm{e}}
 \def\sgn{\mathrm{sgn}}\def\supp{\mathrm{supp}}

\def\rrv{}

\begin{document}

$~$

\bigskip\bigskip

\centerline{\Large\bf Stochastic equations for two-type continuous-state}

\smallskip

\centerline{\Large\bf branching processes in varying environments\footnote{The research is supported by the National Key R\&D Program of China (No. 2020YFA0712901).} }

\bigskip

\centerline{Zenghu Li\footnote{E-mail\,$:$ lizh@bnu.edu.cn. Laboratory of Mathematics and Complex Systems, School of Mathematical Sciences, Beijing Normal University, Beijing 100875, People's Republic of China. }
  \quad
Junyan Zhang$^{*}$\footnote{*corresponding author, E-mail\,$:$ zhangjunyan@mail.bnu.edu.cn. Laboratory of Mathematics and Complex Systems, School of Mathematical Sciences, Beijing Normal University, Beijing 100875, People's Republic of China.} }

\begin{center}
\begin{minipage}{12cm}
 
\begin{center}\textbf{Abstract}\end{center}

\footnotesize

A two-type continuous-state branching process in varying environments is constructed as the pathwise unique solution of a system of stochastic equations driven by time-space noises, where the pathwise uniqueness is derived from a comparison property of solutions. As an application of the main result, we give characterizations of some positive integral functionals of the process in terms of Laplace transforms.

\bigskip



\end{minipage}
\end{center}

\bigskip\bigskip


\section{Introduction}

\setcounter{equation}{0}

Stochastic equations have played important roles in studying continuous-state branching processes (CB-processes) and a number of their generalizations. A general CB-process with immigration (CBI-process) was characterized as the pathwise unique solution to a stochastic equation in Dawson and Li (2006). A flow of discontinuous CB-processes was constructed in Bertoin and Le~Gall (2006) by weak solutions to a stochastic equation; see also Bertoin and Le~Gall (2000, 2003, 2005). The construction was extended to flows of CBI-processes in Dawson and Li (2012) using strong solutions; see also Li (2014) and Li and Ma (2008). For the stable branching CBI-process, a strong stochastic differential equation driven by L\'{e}vy processes was established in Fu and Li (2010). The readers may refer to Bansaye and M\'el\'eard (2015), Li (2020) and Pardoux (2016) for systematic treatments of the theory.

A limit theorem for discrete Galton--Watson branching processes in varying environments was established by Bansaye and Simatos (2015), which leads to CB-processes in varying environments (CBVE-processes). The CBVE-process involves some irregularities at its \textit{bottlenecks}, which are times when the process jumps to zero. In fact, the determination of the behavior of a general CBVE-process at the bottlenecks was left open in Bansaye and Simatos (2015). The problem was settled in Fang and Li (2022) by a construction of the process as the pathwise unique solution to a stochastic integral equation driven by inhomogeneous time-space noises.

Stochastic equations for multi-type CB- and CBI-processes have also been studied. In particular, Watanabe (1969) proved that a two-type continuous CB-process can be constructed as the pathwise unique strong solution of an SDE driven by Brownian motions. Ma (2013) studied the SDE corresponding to a general two-type CBI-process. Barczy et al.\ (2015) constructed multi-type CBI-processes as solutions to vector-valued SDEs.

We are interested in the construction and characterization of the two-type CBVE-process (TCBVE-process). In the recent work Li and Zhang (2025+), we have given a construction of the cumulant semigroup of the TCBVE-process under a first moment condition. The main purpose of this paper is to establish a system of stochastic equations of the process. As an application of the main theorem, we give characterizations of some positive integral functionals of the TCBVE-process in terms of Laplace transforms. The results obtained here will serve as the basis of a further construction of the general TCBVE-process without the moment condition, which will be carried out in a future work.

In the sequel, we denote vectors by boldface letters. Let $\langle\cdot, \cdot\rangle$ and $\|\cdot\|$ denote the Euclidean inner product and norm on $\mathbb{R}^2$, respectively. In the integrals on the real line, we make the convention that, for $t\geq r\in\mathbb{R}$,
\begin{eqnarray*}
	\int_r^t=-\int_t^r=\int_{(r,t]}\text{ and }\int_r^\infty=-\int_\infty^r=\int_{(r,\infty)}.
\end{eqnarray*}

By a \textit{TCBVE-process} (two-type continuous-state branching process in varying environments) we mean an inhomogeneous Markov process $\mbf{X}= \{(X_1(t),X_2(t)): t\ge 0\}$ in the state space $\mbb{R}_+^2$ with transition semigroup $(Q_{r,t})_{t\geq r\ge 0}$ defined by
 \beqlb\label{Laplace transform}
\int_{\mbb{R}_+^2}\e^{-\langle {\bm{\lambda}},\mbf{y}\rangle}Q_{r,t}(\mbf{x},\d\mbf{y})
 =
\e^{-\langle \mbf{x},\mbf{v}_{r,t}({\bm\lambda})\rangle}, \quad \bm\lambda,\mbf{x}\in \mbb{R}_+^2,
 \eeqlb
where $(\mbf{v}_{r,t})_{t\geq r\ge 0}= (v_{1,r,t},v_{2,r,t})_{t\geq r\ge 0}$ is a family of continuous transformations on $\mbb{R}_+^2$.
 \beqlb\label{v-semigroup property}
\mbf{v}_{r,t}= \mbf{v}_{r,s}\circ\mbf{v}_{s,t}, \quad t\ge s\ge r\ge 0,
 \eeqlb
we call $(\mbf{v}_{r,t})_{t\geq r\ge 0}$ the \textit{cumulant semigroup} of $\mbf{X}$.

For $i=1,2$, let $c_i$ be an increasing continuous function on $[0,\infty)$ satisfying $c_i(0)=0$ and let $b_{ii}$ be a c\`adl\`ag function on $[0,\infty)$ satisfying $b_{ii}(0)=0$ and having locally bounded variations. Let $m_i$ be a $\sigma$-finite measure on $(0,\infty)\times (\mbb{R}_+^2\setminus\{\mbf{0}\})$ satisfying
 \beqlb\label{m_i(t)}
m_i(t):= \int_0^t\int_{\mbb{R}_+^2\setminus\{\mbf{0}\}}(z_i^2\mbf 1_{\{\Vert\mbf{z}\Vert\leq 1\}}+z_i\mbf 1_{\{\Vert\mbf{z}\Vert> 1\}}+z_j)m_i(\d s,\d\mbf{z})<\infty, \quad t\geq 0.
 \eeqlb
Moreover, we assume that
\beqlb\label{delta_i}
\delta_i(t)= \Delta b_{ii}(t) + \int_{\mbb{R}_+^2\setminus\{\mbf{0}\}}z_im_i(\{t\},\d\mbf{z})\leq 1, \quad t\geq 0,
\eeqlb
where $\Delta b_{ii}(t)= b_{ii}(t) - b_{ii}(t-)$. For $i,j=1,2$ with $i\neq j$, let $b_{ij}$ be an increasing c\`adl\`ag function on $[0,\infty)$ satisfying $b_{ij}(0)=0$. By Theorem~1.1 in Li and Zhang (2025+), for any $\bm\lambda\in\mbb{R}_+^2$ and $t\geq 0$, there is a unique bounded solution to the following system of backward integral equations
 \beqlb\label{backward eq}
 \begin{split}
 	v_{i,r,t}(\bm\lambda)=\lambda_i &+ \int_r^t v_{j,s,t}(\bm\lambda)b_{ij}(\d s) - \int_r^t v_{i,s,t}(\bm\lambda)b_{ii}(\d s) - \int_r^t v_{i,s,t}(\bm\lambda)^2c_i(\d s) \\
&-\int_r^t\int_{\mbb{R}_+^2\setminus\{\mbf{0}\}}K_i(\mbf{v}_{s,t}(\bm\lambda),\mbf{z})m_i(\d s,\d\mbf{z}), \qquad r\in [0,t],~ i,j=1,2,~ j\neq i,
 \end{split}
 \eeqlb
where
 $$
K_i(\bm\lambda,\mbf{z})=\e^{-\langle \bm\lambda,\mbf{z}\rangle}-1+\lambda_i z_i.
 $$
Moreover, the family $(\mathbf v_{r,t})_{t\geq r\geq 0}$ is the cumulant semigroup of a TCBVE-process.

Suppose that $(\Omega,\mathscr F,\mathscr F_t,\mathbf P)$ is a filtered probability space satisfying the usual hypotheses. For $i=1,2,$ let $W_i(\d s,\d u)$ be a time-space $(\mathscr F_t)$-Gaussian white noise on $(0,\infty)^2$ with intensity $2c_i(\d s)\d u$, and let $M_i(\d s,\d \mathbf z,\d u)$ be a time-space $(\mathscr F_t)$-Poisson random measure on $(0,\infty)\times (\mathbb R_+^2\backslash\{\mathbf 0\})\times(0,\infty)$ with intensity $m_i(\d s,\d \mathbf z)\d u$. Suppose that $W_1(\d s,\d u)$, $W_2(\d s,\d u)$, $M_1(\d s,\d \mathbf z,\d u)$ and $M_2(\d s,\d \mathbf z,\d u)$ are independent of each other. We denote by $\tilde M_i(\d s,\d \mathbf z,\d u)$ the compensated measure of $M_i(\d s,\d \mathbf z,\d u)$. Let us consider the following system of stochastic equations
 \beqlb\label{SIE}
 \begin{split}
X_i(t)&=X_i(0)+\int_0^t\int_0^{X_i(s-)}W_i(\d s,\d u)-\int_0^t X_i(s-)b_{ii}(\d s)+\int_0^t X_j(s-) b_{ji}(\d s) \\
 &\quad\
+\int_0^t\int_{\mathbb R_+^2\backslash\{\mathbf 0\}}\int_0^{X_i(s-)}z_i\tilde M_i(\d s,\d \mathbf z,\d u)+\int_0^t\int_{\mathbb R_+^2\backslash\{\mathbf 0\}}\int_0^{X_j(s-)}z_i M_j(\d s,\d \mathbf z,\d u),
 \end{split}
 \eeqlb
where $r\in [0,t],~ i,j=1,2,~ j\neq i$. The main results of the paper are the following:

\begin{theorem}\label{SIE solution} There is a pathwise unique solution $\{\mathbf X(t):t\geq 0\}$ to \eqref{SIE} and the solution is a TCBVE-process with transition semigroup $(Q_{r,t})_{t\geq r}$ defined by \eqref{Laplace transform} and \eqref{backward eq}. \end{theorem}

\begin{theorem}\label{positive integral functionals} Let $\mathbf X=(\Omega,\mathscr F,\mathscr F_{r,t},\mathbf X(t),\mathbf P_{r,\mathbf x})$ be a TCBVE-process started from time $r\ge 0$. For $i=1,2$, let $\zeta_i$ be a $\sigma$-finite measure on $\mathbb R_+$ satisfying $\zeta_i(B)<\infty$ for every bounded Borel set $B\subset \mathbb R_+$. Then we have
 \begin{equation*}
\mathbf P_{r,\mathbf x}\exp\bigg\{-\sum_{i=1}^2\int_{[r,t]} X_i(s)\zeta_i(\d s)\bigg\}=\exp\{-\langle\mathbf x,\mathbf w_{r,t}\rangle\}, \quad t\geq r,\, \mbf x\in \mathbb R_+^2,
 \end{equation*}
where $r\mapsto\mathbf w_{r,t}$ is the unique positive solution to
 \beqnn
w_{i,r,t}=\int_{[r,t]}\zeta_i(\d s)\ar-\ar\int_r^t w_{i,s,t}b_{ii}(\d s)+\int_r^t w_{j,s,t}b_{ij}(\d s)-\int_r^t w_{i,s,t}^2 c_{i}(\d s)\\
 \ar-\ar
\int_r^t\int_{\mathbb R_+^2\backslash\{\mathbf 0\}}K_i(\mathbf w_{s,t},\mathbf z)m_i(\d s,\d \mathbf z),\qquad r\in [0,t],~ i,j=1,2,~ j\neq i.
	\eeqnn
\end{theorem}

The rest of the paper is organized as follows. In Section 2, we establish a comparison property for the solutions to \eqref{SIE}, which implies the pathwise uniqueness of the solution with arbitrary initial state. In Section 3, we prove that any TCBVE-process with transition semigroup $(Q_{r,t})_{t\geq r}$ defined by \eqref{Laplace transform} and \eqref{backward eq} is a weak solution to \eqref{SIE}. Based on those results, the proofs of Theorems~\ref{SIE solution} and~\ref{positive integral functionals} are given in Section 4 and Section 5, respectively.

\section{Comparison property}

\setcounter{equation}{0}

In this section, we prove a comparison property of solutions to the stochastic equation system \eqref{SIE}, which implies the pathwise uniqueness for the system. For a c\`adl\`ag function $f$, let $J_f:=\{s>0:|\triangle f(s)|>0\}$ denote the discontinuous points of $f$, which is at most a countable set. Then for $i=1,2$, we can define two $(\mathscr F_t)$-Poisson random measures $M_{i,c}(\d s,\d \mathbf z,\d u):=\mathbf 1_{J_{m_i}^c}(s)M_i(\d s,\d \mathbf z,\d u)$ and $M_{i,d}(\d s,\d \mathbf z,\d u):=\mathbf 1_{J_{m_i}}(s)M_i(\d s,\d \mathbf z,\d u)$ with intensities $m_{i,c}(\d s,\d \mathbf z)\d u$ and $m_{i,d}(\d s,\d \mathbf z)\d u$, respectively. Moreover, $M_{i,c}$ and $M_{i,d}$ are independent of each other as they have disjoint supports. In this setting, we can rewrite \eqref{SIE} into
\beqlb\label{4 SIE}
	\begin{split}
		X_i&(t)=X_i(0)+\int_0^t\int_0^{X_i(s-)}W_i(\d s,\d u)-\int_0^t X_i(s-)b_{ii}(\d s)+\int_0^t X_j(s-)\bar b_{ji}(\d s)\\
		&+\int_0^t\int_{\mathbb R_+^2\backslash\{\mathbf 0\}}\int_0^{X_i(s-)}z_i\tilde M_{i,c}(\d s,\d \mathbf z,\d u)
		+\int_0^t\int_{\mathbb R_+^2\backslash\{\mathbf 0\}}\int_0^{X_i(s-)}z_i\tilde M_{i,d}(\d s,\d \mathbf z,\d u)\\
		&+\int_0^t\int_{\mathbb R_+^2\backslash\{\mathbf 0\}}\int_0^{X_j(s-)}z_i\tilde M_{j,c}(\d s,\d \mathbf z,\d u)
		+\int_0^t\int_{\mathbb R_+^2\backslash\{\mathbf 0\}}\int_0^{X_j(s-)}z_i\tilde M_{j,d}(\d s,\d \mathbf z,\d u),
	\end{split}
\eeqlb
where
\beqlb\label{b_ij-}
\bar{b}_{ij}(t)=b_{ij}(t)+\int_0^t\int_{\mbb{R}_+^2\setminus\{\mbf{0}\}}z_jm_i(\d s,\d\mbf{z}), \quad t\ge 0.
\eeqlb

\bproposition\label{EX_(t) estimate} Let $\{\mathbf X(t):t\geq 0\}$ be a solution to \eqref{4 SIE} satisfying $\mathbf E[X_i(0)]<\infty$. Let $\tau_{i,k}:=\inf\{t\geq 0:X_i(t)\geq k\}$ for $k\geq 1,~ i=1,2$. Then $\tau_{i,k}\to\infty$ a.s. as $k\to\infty$. Moreover, for $t\geq 0,~ i=1,2$ and $k\geq 1$ we have
 \begin{equation}
k\mathbf P\{\tau_{i,k}\leq t\}\leq\mathbf E[X_i(t)] \leq \mathbf E\mathbf [X_i(0)]\exp\{\alpha(t)\}+\mathbf E[X_j(0)]\bar b_{ji}(t)\exp\{\Vert b_{jj}\Vert(t)+\alpha(t)\},
 \end{equation}
where
 \begin{equation*}
\alpha(t)=\bar b_{12}(t)\bar b_{21}(t)\exp\{\Vert b_{jj}\Vert (t)\}+\Vert b_{ii}\Vert(t).
 \end{equation*}
\eproposition

\bproof It is obvious by Fubini's Theorem and Gronwall's inequality. \eproof

\blemma\label{disjoint jumps} Suppose that $N_1(\d s,\d\mathbf z)$ and $N_2(\d s,\d\mathbf z)$ are independent Poisson random measure on $(0,\infty)\times (\mathbb R^n\backslash\{\mbf 0\})$ with intensities $\d s n_1(\d\mathbf z)$ and $\d s n_2(\d \mathbf z)$, respectively. Suppose that $f_1,f_2$ are Borel functions on $\mathbb R^n$. Let
 \beqlb
N_i(t)=\int_0^t\int_{\mathbb R^n\backslash\{\mbf 0\}}f_i(\mathbf z)N_i(\d s,\d\mathbf z),\quad i=1,2.
 \eeqlb
Let $J_{N_i}=\{t> 0:|\triangle N_i(t)|>0\}$. Then $\mathbf P(J_{N_1}\cap J_{N_2}\neq\emptyset)=0$. \elemma

\bproof Suppose that $\mathscr G_t$ is the natural filtration generated by $N_1$ satisfying usual hypotheses. Let $\mathscr G_\infty=\sigma(\cup_{t\geq 0} \mathscr G_t)$. Fix $\omega\in\Omega$, it is obvious that $J_{N_1}(\omega)$ is at most a countable set. Then
 \beqnn
\mathbf P(J_{N_1}\cap J_{N_2}\neq\emptyset)
 \ar=\ar
\mathbf P(\mathbf P(J_{N_1}\cap J_{N_2}\neq\emptyset|\mathscr G_\infty))\\
\ar=\ar\mathbf P(\mathbf P(A\cap J_{N_2}\neq\emptyset|\mathscr G_\infty)|_{A=J_{N_1}\in \mathscr G_\infty})\\
 \ar=\ar
\mathbf P(\mathbf P(A\cap J_{N_2}\neq\emptyset)|_{A=J_{N_1}\in \mathscr G_\infty}) =0,
 \eeqnn
where we use $N_1$ and $N_2$ are independent for the third equality. \eproof

Now we give the comparison property. For simplicity, we only consider the case under an additional second moment integrability condition.

\bproposition\label{comparison property} Assume that
 \begin{equation}\label{additional integrability condition}
\int_0^t\int_{\mathbb R_+^2\backslash\{\mathbf 0\}}(z_1^2+z_2^2)(m_1+m_2)(\d s,\d \mathbf z)<\infty,\quad t\geq 0.
 \end{equation}
If $\{\mathbf X^{(1)}(t):t\geq 0\}$ and $\{\mathbf X^{(2)}(t):t\geq 0\}$ are two solutions to \eqref{4 SIE} satisfying $\mathbf P\{X_i^{(1)}(0)\leq X_i^{(2)}(0),i=1,2\}=1$, then we have $\mathbf P\{X_i^{(1)}(t)\leq X_i^{(2)}(t)\text{ for every }t\geq 0,i=1,2\}=1$, which implyies the pathwise uniqueness of solution holds for \eqref{4 SIE}. \eproposition

\bproof By passing to the conditional law $\mathbf P(\cdot|\mathscr F_0)$ if it is necessary, we may assume both $\mathbf X^{(1)}(0)$ and $\mathbf X^{(2)}(0)$ are deterministic. Then the results of Proposition \ref{EX_(t) estimate} are valid. For each integer $n\geq 0$, define $a_n=\exp\{-n(n+1)/2\}$. Note that $a_n\to 0$ decreasingly as $n\to\infty$ and
 \begin{equation*}
\int_{a_n}^{a_{n-1}}2(nz)^{-1}\d z=2n^{-1}\log\left(\frac{a_{n-1}}{a_n}\right)=2.
 \end{equation*}
Then there is a positive continuous function $x\mapsto g_n(x)$ supported by $(a_n,a_{n-1})$ so that $$\int_{a_n}^{a_{n-1}} g_n(x)\d x=1$$ and $g_n(x)\leq 2(nx)^{-1}$ for every $x>0$. For $n\geq 0$ and $z\in\mathbb R$, let
 \begin{equation*}
f_n(z)=\int_0^{z\vee 0}\d y\int_0^y g_n(x)\d x.
 \end{equation*}
It is easy to see that $f_n(z)=0$ for $z< 0$ and $f_n(z)\to z$ increasingly for $z\geq 0$ as $n\to\infty$. For simplicity, let $(i,j)=(1,2)$ or $(2,1)$. Setting $$Y_i(t)=X_i^{(1)}(t)-X_i^{(2)}(t),$$from \eqref{4 SIE} we have
 \begin{equation}\label{Y_i(t)}
 \begin{split}
Y_i(t)&=Y_i(0)+\int_0^t\int_{\underline X_i(s-)}^{\overline X_i(s-)} q_i(s-)W_i(\d s,\d u)-\int_0^t Y_i(s-)b_{ii}(\d s)\\
 &\quad
+\int_0^t Y_j(s-)\bar b_{ji}(\d s)+\int_0^t\int_{\mathbb R_+^2\backslash\{\mathbf 0\}}\int_{\underline X_i(s-)}^{\overline X_i(s-)}z_iq_i(s-)(\tilde M_{i,c}+\tilde M_{i,d})(\d s,\d \mathbf z,\d u)\\
 &\quad
+\int_0^t\int_{\mathbb R_+^2\backslash\{\mathbf 0\}}\int_{\underline X_j(s-)}^{\overline X_j(s-)}z_iq_j(s-)(\tilde M_{j,c}+\tilde M_{j,d})(\d s,\d \mathbf z,\d u),
 \end{split}
 \end{equation}
where $$\underline X_i(s-)=X_i^{(1)}(s-)\wedge X_i^{(2)}(s-),\ \overline X_i(s-)=X_i^{(1)}(s-)\vee X_i^{(2)}(s-)$$ and $q_i(s)=\mathbf 1_{\{Y_i(s-)>0\}}-\mathbf 1_{\{Y_i(s-)<0\}}$. Let $$D_z f_n(x)=f_n(x+z)-f_n(x)-f_n'(x)z.$$ By It\^o's formula,
	\begin{equation}\label{f_n(Y_i(t))}
		\begin{split}
			f_n(Y_i(t))&=f_n(Y_i(0))+\int_0^t f_n'(Y_i(s-))\d Y_i(s)\\
			&\quad+\frac{1}{2}\int_0^t f_n''(Y_i(s-))\d [Y_i,Y_i]^{c}(s)+\sum_{s\in(0,t]}D_{\triangle Y_i(s)}f_n(Y_i(s-))\\
			&=-\int_0^t f_n'(Y_i(s-)) Y_i(s-)b_{ii}(\d s)+\int_0^t f_n'(Y_i(s-)) Y_j(s-)\bar b_{ji}(\d s)\\
			&\quad+\int_0^t|Y_i(s-)|f_n''(Y_i(s-))c_i(\d s)+\sum_{s\in(0,t]}D_{\triangle Y_i(s)}f_n(Y_i(s-))+M_i(t),
		\end{split}
	\end{equation}	
	where
	\begin{eqnarray*}
		M_i(t)\ar=\ar\int_0^t\int_{\underline X_i(s-)}^{\overline X_i(s-)} f_n'(Y_i(s-))q_i(s-)W_i(\d s,\d u)\\
			\ar\ar+\int_0^t\int_{\mathbb R_+^2\backslash\{\mathbf 0\}}\int_{\underline X_i(s-)}^{\overline X_i(s-)}f_n'(Y_i(s-))z_iq_i(s-)(\tilde M_{i,c}+\tilde M_{i,d})(\d s,\d \mathbf z,\d u)\\
			\ar\ar+\int_0^t\int_{\mathbb R_+^2\backslash\{\mathbf 0\}}\int_{\underline X_j(s-)}^{\overline X_j(s-)}f_n'(Y_i(s-))z_iq_j(s-)(\tilde M_{j,c}+\tilde M_{j,d})(\d s,\d \mathbf z,\d u).
	\end{eqnarray*}
	We can easily obtain that $\{M_i(t):t\geq 0\}$ is a square-integrable martingale by Proposition \ref{EX_(t) estimate}. In fact,
	\begin{eqnarray*}
		\mathbf E\left[M_i(t)^2\right]\ar\leq\ar 8\mathbf E\left[\int_0^t|Y_i(s-)|c_i(\d s)\right]+4\mathbf E\left[\int_0^t\int_{\mathbb R_+^2\backslash\{\mathbf 0\}}z_i^2|Y_i(s-)|m_i(\d s,\d \mathbf z)\right]\\
			\ar\ar+4\mathbf E\left[\int_0^t\int_{\mathbb R_+^2\backslash\{\mathbf 0\}}z_i^2|Y_j(s-)|m_j(\d s,\d \mathbf z)\right].
	\end{eqnarray*}
	Let $\supp(M_{i,c})\subset(0,\infty)\times(\mathbb R_+^2\backslash\{\mathbf 0\})\times(0,\infty)$ denote the countable support of Poisson random measure $M_{i,c}(\d s,\d \mathbf z,\d u)$. Let $J_{M_{i,c}}$ denote the projection of $\supp(M_{i,c})$ to the temporal axis. Denote the set of jump times of $Y_i(t)$ by $J_{Y_i}$ and $B_i:=J_{b_{ji}}\cup J_{b_{ii}}\cup J_{m_1}\cup J_{m_2}$. Then by \eqref{Y_i(t)} and Lemma \ref{disjoint jumps}, we can divide $J_{Y_i}$ into three disjoint sets a.s., that is,
	\beqlb
		J_{Y_i}=B_i\bigcup J_{M_{i,c}}\bigcup J_{M_{j,c}}.
	\eeqlb
	For $s\in J_{M_{k,c}},~ k=1,2$, the process $\{Y_i(t):t\geq 0\}$ makes a jump of size
	\beqnn
	\triangle Y_i(s)=\sgn(Y_k(s-))z_i\mathbf 1_{\{\underline X_k(s-)<u\leq \overline X_k(s-)\}},
	\eeqnn
	where $\mathbf z\in\mathbb R_+^2\backslash\{\mathbf 0\}$, $u>0$ satisfying $(s,\mathbf z,u)\in\supp(M_{k,c})$, so $\sgn(\triangle Y_i(s))=\sgn(Y_k(s-))$. At time $s\in B_i$, the process $\{Y_i(t):t\geq 0\}$ jumps by
	\beqlb
		\begin{aligned}
			\triangle Y_i(s)&=-Y_i(s-)\triangle b_{ii}(s)+\int_{\mathbb R_+^2\backslash\{\mathbf 0\}}\int_{\underline X_i(s-)}^{\overline X_i(s-)} z_i q_i(s-)\tilde M_{i,d}(\{s\},\d \mathbf z,\d u)\\
			&\quad+Y_j(s-)\triangle\overline b_{ji}(s)+\int_{\mathbb R_+^2\backslash\{\mathbf 0\}}\int_{\underline X_j(s-)}^{\overline X_j(s-)} z_i q_j(s-)\tilde M_{j,d}(\{s\},\d \mathbf z,\d u)\\
			&=-Y_i(s-)\delta_i(s)+\int_{\mathbb R_+^2\backslash\{\mathbf 0\}}\int_{\underline X_i(s-)}^{\overline X_i(s-)} z_i q_i(s-) M_{i,d}(\{s\},\d \mathbf z,\d u)\\
			&\quad+Y_j(s-)\triangle b_{ji}(s)+\int_{\mathbb R_+^2\backslash\{\mathbf 0\}}\int_{\underline X_j(s-)}^{\overline X_j(s-)} z_i q_j(s-) M_{j,d}(\{s\},\d \mathbf z,\d u),
		\end{aligned}
	\eeqlb
	where $\delta_i$ is defined by \eqref{delta_i}. Together with the property of $f_n$, we can rewrite \eqref{f_n(Y_i(t))} into
	\beqlb\label{f_n(Y_i(t))rewrite}
		f_n(Y_i(t))&\!\!\!=\!\!\!&M_i(t)-\int_0^t f_n'(Y_i(s-)) Y_i(s-)b_{ii}(\d s)
		+\int_0^t|Y_i(s-)|f_n''(Y_i(s-))c_i(\d s)\notag\\
		&&\quad+\int_0^t f_n'(Y_i(s-)) Y_j(s-) b_{ji}(\d s)
		+\sum_{s\in(0,t]}D_{\triangle Y_i(s)}f_n(Y_i(s-))\mathbf 1_{B_i}(s)\notag\\
		&&\quad+\int_0^t\int_{\mathbb R_+^2\backslash\{\mathbf 0\}}\int_{\underline X_j(s-)}^{\overline X_j(s-)} z_iq_j(s-)f_n'(Y_i(s-)) m_{j,c}(\d s,\d \mathbf z)\d u\notag\\
		&&\quad+\int_0^t\int_{\mathbb R_+^2\backslash\{\mathbf 0\}}\int_{\underline X_j(s-)}^{\overline X_j(s-)} z_iq_j(s-)f_n'(Y_i(s-)) m_{j,d}(\d s,\d \mathbf z)\d u\\
		&&\quad+\int_0^t\int_{\mathbb R_+^2\backslash\{\mathbf 0\}}\int_{\underline X_i(s-)}^{\overline X_i(s-)} D_{z_i}f_n(Y_i(s-))\mathbf 1_{\{Y_i(s-)>0\}} M_{i,c}(\d s,\d \mathbf z,\d u)\notag\\
		&&\quad+\int_0^t\int_{\mathbb R_+^2\backslash\{\mathbf 0\}}\int_{\underline X_j(s-)}^{\overline X_j(s-)} D_{q_j(s-) z_i}f_n(Y_i(s-))M_{j,c}(\d s,\d \mathbf z,\d u)\notag\\
		&\!\!\!=:\!\!\!&M_i(t)+\sum_{l=1}^8 I_l,\notag
	\eeqlb
	where $I_l,1\leq l\leq 8$ is defined by the $(l+1)$th term of the right hand side of \eqref{f_n(Y_i(t))rewrite}. By elementary calculations, we have
	\begin{eqnarray*}
		\mathbf E[I_1]\ar\leq\ar\int_0^t\mathbf E[Y_i(s-)\vee 0]\Vert b_{ii}\Vert(\d s).
	\end{eqnarray*}
	Denote $\triangle_z f_n(x)=f_n(x+z)-f_n(x)$, we can easily see $\triangle_z f_n(x)\leq 0$ when $z\leq 0$. Then we have
	\begin{eqnarray*}
		\mathbf E[I_5+I_8]\ar=\ar\mathbf E\left[\int_0^t\int_{\mathbb R_+^2\backslash\{\mathbf 0\}}\triangle_{q_j(s-)z_i}f_n(Y_i(s-))|Y_j(s-)|m_{j,c}(\d s,\d \mathbf z) \right]\\
		\ar\leq\ar \mathbf E\left[\int_0^t\int_{\mathbb R_+^2\backslash\{\mathbf 0\}}\triangle_{z_i}f_n(Y_i(s-))|Y_j(s-)|\mathbf 1_{\{Y_j(s-)\geq 0\}} m_{j,c}(\d s,\d \mathbf z) \right]\\
		\ar\leq\ar \int_0^t\int_{\mathbb R_+^2\backslash\{\mathbf 0\}}\mathbf E[Y_j(s-)\vee 0]z_i m_{j,c}(\d s,\d \mathbf z).
	\end{eqnarray*}
	As in the proof of Theorem 8.2 in Li (2020) one sees $|xf_n''(x)|\leq 2/n$ and $|xD_z f_n(x)|\leq z^2/n$ for $x,z\geq 0$. Then we have $\mathbf E[I_2]\leq 2c_i(t)/n\to 0$ as $n\to\infty$ and
	\begin{eqnarray*}
		\mathbf E[I_7]\ar=\ar\mathbf E\left[\int_0^t\int_{\mathbb R_+^2\backslash\{\mathbf 0\}} Y_i(s-)D_{z_i}f_n(Y_i(s-))\mathbf 1_{\{Y_i(s-)>0\}} m_{i,c}(\d s,\d \mathbf z) \right]\\
		\ar\leq\ar \frac{1}{n}\int_0^t\int_{\mathbb R_+^2\backslash\{\mathbf 0\}} z_i^2 m_{i,c}(\d s,\d \mathbf z)\to 0,\quad n\to\infty.
	\end{eqnarray*}
	Denote
	\begin{eqnarray*}
		I_{4,1}\ar=\ar\sum_{s\in(0,t]}D_{\triangle Y_i(s)}f_n(Y_i(s-))\mathbf 1_{B_i}(s)\mathbf 1_{\{Y_j(s-)\leq 0\}},\\
		I_{4,2}\ar=\ar\sum_{s\in(0,t]}D_{\triangle Y_i(s)}f_n(Y_i(s-))\mathbf 1_{B_i}(s)\mathbf 1_{\{Y_j(s-)>0\}}
	\end{eqnarray*}
	and
	\beqlb\label{triangle Y_i,1(s)}
	\begin{aligned}
		\triangle Y_{i,1}(s)&=-Y_i(s-)\delta_i(s)+\int_{\mathbb R_+^2\backslash\{\mathbf 0\}}\int_{\underline X_i(s-)}^{\overline X_i(s-)} z_i q_i(s-) M_{i,d}(\{s\},\d \mathbf z,\d u),\\
		\triangle Y_{i,2}(s)&=Y_j(s-)\triangle b_{ji}(s)+\int_{\mathbb R_+^2\backslash\{\mathbf 0\}}\int_{\underline X_j(s-)}^{\overline X_j(s-)} z_i q_j(s-) M_{j,d}(\{s\},\d \mathbf z,\d u),
	\end{aligned}
	\eeqlb
	so $\sgn(Y_i(s-)+\triangle Y_{i,1}(s))=\sgn(Y_i(s-))$, implying $D_{\triangle Y_{i,1}(s)}f_n(Y_i(s-))\mathbf 1_{\{Y_i(s-)\leq 0\}}=0$.
	For $z_2\leq 0$, it's easy to see $D_{z_1+z_2}f_n(x)\leq D_{z_1}f_n(x)-f_n'(x)z_2$. For $s\in B_i$, we have
	\begin{equation*}
		D_{\triangle Y_i(s)}f_n(Y_i(s-))\mathbf 1_{\{Y_j(s-)\leq 0\}}\leq [D_{\triangle Y_{i,1}(s)}f_n(Y_i(s-))-f_n'(Y_i(s-))\triangle Y_{i,2}(s)]\mathbf 1_{\{Y_j(s-)\leq 0\}},
	\end{equation*}
	then
	\begin{eqnarray*}
		\mathbf E[I_{4,1}+I_3+I_6]\ar\leq\ar \mathbf E\left[\sum_{s\in(0,t]}D_{\triangle Y_{i,1}(s)}f_n(Y_i(s-))\mathbf 1_{\{Y_j(s-)\leq 0\}}\mathbf 1_{B_i}(s) \right]\\
		\ar\ar+\mathbf E\left[\int_0^t f_n'(Y_i(s-))Y_j(s-)\mathbf 1_{\{Y_j(s-)>0\}} b_{ji}(\d s) \right]\\
		\ar\ar+\mathbf E\left[\int_0^t\int_{\mathbb R_+^2\backslash\{\mathbf 0\}}\int_{\underline X_j(s-)}^{\overline X_j(s-)} z_i q_j(s-)\mathbf 1_{\{Y_j(s-)>0\}} m_{j,d}(\d s,\d \mathbf z)\d u \right]\\
		\ar\leq\ar \int_0^t \mathbf E[Y_j(s-)\vee 0] b_{ji}(\d s)+\int_0^t\int_{\mathbb R_+^2\backslash\{\mathbf 0\}} \mathbf E[Y_j(s-)\vee 0]z_i m_{j,d}(\d s,\d \mathbf z)\\
		\ar\ar +\mathbf E\left[\sum_{s\in(0,t]}D_{\triangle Y_{i,1}(s)}f_n(Y_i(s-))\mathbf 1_{\{Y_j(s-)\leq 0\}}\mathbf 1_{B_i}(s) \right].
	\end{eqnarray*}
	For $z_2>0$, it's easy to see $D_{z_1+z_2}f_n(x)\leq D_{z_1}f_n(x)+z_2$. For $s\in B_i$, we have
	\begin{equation*}
		D_{\triangle Y_i(s)}f_n(Y_i(s-))\mathbf 1_{\{Y_j(s-)>0\}}\leq [D_{\triangle Y_{i,1}(s)}f_n(Y_i(s-))+\triangle Y_{i,2}(s)]\mathbf 1_{\{Y_j(s-)>0\}},
	\end{equation*}
	then
	\begin{eqnarray*}
		\mathbf E[I_{4,2}]\ar\leq\ar \int_0^t \mathbf E[Y_j(s-)\vee 0] b_{ji}(\d s)+\int_0^t\int_{\mathbb R_+^2\backslash\{\mathbf 0\}} \mathbf E[Y_j(s-)\vee 0]z_i m_{j,d}(\d s,\d \mathbf z)\\
		\ar\ar+\mathbf E\left[\sum_{s\in(0,t]}D_{\triangle Y_{i,1}(s)}f_n(Y_i(s-))\mathbf 1_{\{Y_j(s-)> 0\}}\mathbf 1_{B_i}(s) \right].
	\end{eqnarray*}
	By \eqref{triangle Y_i,1(s)}, for $s\in B_i$ we have $\triangle Y_{i,1}(s)=g_{s,i}(X^{(1)}(s-))-g_{s,i}(X^{(2)}(s-))$, where
	\begin{equation*}
		g_{s,i}(x_1,x_2)=\int_{\mathbb R_+^2\backslash\{\mathbf 0\}}\int_0^{x_i}z_iM_{i,d}(\{s\},\d \mathbf z,\d u)-x_i\delta_i(s).
	\end{equation*}
	Since $x_i\mapsto x_i+g_{s,i}(x_1,x_2)$ is an increasing function, by Lemma 3.1 in Li and Pu (2012), we have, for $s\in B_i$,
	\begin{eqnarray*}
		|D_{\triangle Y_{i,1}(s)} f_n(Y_i(s-))|\ar\leq\ar 2n^{-1}|Y_i(s-)|^{-1}[|\triangle Z_i(s)|+|Y_i(s-)||\triangle b_{ii}(s)|]^2\\
		\ar\leq\ar 4n^{-1}[|Y_i(s-)|^{-1}|\triangle Z_i(s)|^2+|Y_i(s-)||\triangle b_{ii}(s)|^2],
	\end{eqnarray*}
	where
	\begin{equation*}
		\triangle Z_i(s)=\int_{\mathbb R_+^2\backslash\{\mathbf 0\}}\int_{\underline X_i(s-)}^{\overline X_i(s-)} z_i\tilde M_{i,d}(\{s\},\d \mathbf z,\d u),\quad i=1,2.
	\end{equation*}
	Taking expectation, we have
	{\small
	\begin{eqnarray*}
		\ar\ar\!\!\!\!\!\!\!\!\!\mathbf E\left[\sum_{s\in(0,t]}D_{\triangle Y_{i,1}(s)}f_n(Y_i(s-))\mathbf 1_{B_i}(s) \right]\leq 4n^{-1}\mathbf E\left[\sum_{s\in(0,t]}Y_i(s-)^{-1}|\triangle Z_i(s)|^2\mathbf 1_{\{Y_i(s-)> 0\}} \mathbf 1_{J_{m_i}}(s) \right]\\
		\ar\ar\quad+4n^{-1}\mathbf E\left[\sum_{s\in(0,t]}Y_i(s-)|\triangle b_{ii}(s)|^2 1_{\{Y_i(s-)> 0\}} \mathbf 1_{J_{b_{ii}}}(s) \right] \\
		\ar\ar\leq 4n^{-1}\int_0^t\int_{\mathbb R_+^2\backslash\{\mathbf 0\}}z_i^2m_{i,d}(\d s,\d \mathbf z)+4n^{-1}\sum_{s\in(0,t]}\mathbf E[Y_i(s-)\vee 0]|\triangle b_{ii}(s)|^2\to 0,\ n\to\infty.\\		
	\end{eqnarray*}
	}Let $n\to\infty$ in \eqref{f_n(Y_i(t))rewrite}, by Fatou's Lemma, we have
	\begin{eqnarray*}
		\mathbf E[Y_i(t)\vee 0]\ar\leq\ar \int_0^t \mathbf E[Y_i(s-)\vee 0]\Vert b_{ii}\Vert (\d s)+2\int_0^t \mathbf E[Y_j(s-)\vee 0] b_{ji}(\d s)\\
		\ar\ar+3\int_0^t\int_{\mathbb R_+^2\backslash\{\mathbf 0\}} \mathbf E[Y_j(s-)\vee 0] z_im_j(\d s,\d \mathbf z).
	\end{eqnarray*}
Then by Gronwall's inequality, we have $\mathbf E[Y_i(t)\vee 0]=0$ for every $t\geq 0$, implying the desired comparison property. \eproof

\section{Weak solution}

\setcounter{equation}{0}

In this section, we prove that any TCBVE-process with transition semigroup defined by \eqref{Laplace transform} and \eqref{backward eq} is a weak solution to the stochastic equation system \eqref{SIE}.

\bproposition\label{composite Levy-Ito representation} Let $(\Omega,\mathscr F,\mathbf P)$ be a complete probability space with the sub-$\sigma$-algebra $\mathscr G\subset\mathscr F$. Suppose that $\bm\beta\in\mathbb R^2$ and $(z_1\wedge z_1^2+z_2\wedge z_2^2)\gamma(\d \mathbf z)$ is a finite measure on $\mathbb R_+^2\backslash\{\mathbf 0\}$. Let $\Psi$ be a function on $\mathbb R_+^2$ given by
 \begin{equation*}
\Psi(\bm\lambda)= \langle\bm\beta,\bm\lambda\rangle + \int_{\mathbb R_+^2\backslash\{\mathbf 0\}} K(\bm\lambda,\mathbf z)\gamma(\d \mathbf z), \quad \lambda\in\mathbb R_+^2,
 \end{equation*}
where
 $$
K(\bm\lambda,\mathbf z)=e^{-\langle\bm\lambda,\mathbf z\rangle}-1+\langle\bm\lambda,\mathbf z\rangle.
 $$
Suppose that $(\xi,\mathbf Y)$ is a random vector taking values in $[0,\infty)\times\mathbb R^2$ such that $\xi$ is $\mathscr G$-measurable and
\begin{equation}
	\mathbf E(e^{-\langle\bm\lambda,\mathbf Y\rangle}|\mathscr G)=\exp\{\xi\Psi(\bm\lambda)\},\quad \bm\lambda\in\mathbb R_+^2.
\end{equation}
Then on an extension of the probability space there exists a Poisson random measure $N(\d \mathbf z,\d u)$ on $\big(\mathbb R_+^2\backslash\{\mathbf 0\}\big)\times(0,\infty)$ with intensity $\gamma(\d \mathbf z)\d u$ such that $N$ is independent of $\mathscr G$ and a.s.
\begin{equation}
	\mathbf Y=-\bm\beta\xi+\int_{\mathbb R_+^2\backslash\{\mathbf 0\}}\int_0^\xi \mathbf z\tilde N(\d \mathbf z,\d u).
\end{equation}
\eproposition

\bproposition\label{TCBVE realization} The TCBVE-process with transition semigroup $(Q_{r,t})_{t\geq r}$ has a c\`adl\`ag semimartingale realization $\{(\mathbf X(t),\mathscr F_t):t\geq 0\}$ with the filtration satisfying the usual hypotheses. For such a realization and $t\in D:= J_{\bar b_{12}}\cup J_{\bar b_{21}}\cup J_{b_{11}}\cup J_{b_{22}}\cup J_{m_1}\cup J_{m_2}$, we have
 \begin{equation}
\mathbf E(e^{-\langle\bm\lambda,\triangle \mathbf X(t)\rangle}|\mathscr F_{t-})=e^{\langle\bm\lambda-\mathbf v_{t-,t}(\bm\lambda),\mathbf X(t-)\rangle},\quad \bm\lambda\in\mathbb R_+^2,
 \end{equation}
where $\triangle\mathbf X(t)=\mathbf X(t)-\mathbf X(t-)$ and
 \begin{equation}
\lambda_i-v_{i,t-,t}(\bm\lambda)=\lambda_i\triangle b_{ii}(t)-\lambda_j\triangle\overline b_{ij}(t)+\int_{\mathbb R_+^2\backslash\{\mathbf 0\}}K(\bm\lambda,\mathbf z)m_i(\{t\},\d \mathbf z), \quad i,j=1,2,~i\neq j.
 \end{equation}
\eproposition

We omit the proofs of the above propositions since they are quite similar to those in the one-dimensional case given in Fang and Li (2022).

\bproposition\label{SIE's weak solution} Suppose that \eqref{additional integrability condition} holds. Then the TCBVE-process with transition semigroup $(Q_{r,t})_{t\geq r}$ defined by \eqref{Laplace transform} and \eqref{backward eq} is a weak solution to \eqref{4 SIE}. \eproposition

\bproof By Proposition \ref{TCBVE realization}, there is a c\`adl\`ag semimartingale realization $\{(\mathbf X(t),\mathscr F_t):t\geq 0\}$ of the TCBVE-process with conservative cumulant semigroup $(\mathbf v_{r,t}(\bm\lambda))_{t\geq r}$ defined by \eqref{backward eq}. Let $N_0(\d s,\d \mathbf z)$ be the optional random measure on $(0,\infty)\times\mathbb R^2$ by
	\begin{equation*}
		N_0(\d s,\d \mathbf z):=\sum_{s>0}\mathbf 1_{\{\triangle\mathbf X(s)\neq \mathbf 0\}}\delta_{(s,\triangle\mathbf X(s))}(\d s,\d \mathbf z).
	\end{equation*}
	Denote the predictable compensator of $N_0(\d s,\d \mathbf z)$ by $\hat N_0(\d s,\d \mathbf z)$ and let $\tilde N_0(\d s,\d \mathbf z)=N_0(\d s,\d \mathbf z)-\hat N_0(\d s,\d \mathbf z)$ be the compensated measure; see, e.g., Dellacherie and Meyer (1982), page 370-372. Then for $i=1,2$,  we have the decomposition
	\begin{equation}\label{X_i(t)decomposition}
		X_i(t)=X_i(0)+M_i(t)+A_i(t)+\int_0^t\int_{\mathbb R^2} z_i\tilde N_0(\d s,\d \mathbf z),
	\end{equation}
	where $\{M_i(t):t\geq 0\}$ is a continuous local martingale satisfying $M_i(0)=0$ and $\{A_i(t):t\geq 0\}$ is a process with bounded variation satisfying $A_i(0)=0$. Let $\{C_i(t):t\geq 0\}$ be the quadratic variation process of $M_i(t)$. For $\mathbf x,\bm\lambda\in\mathbb R_+^2$, define a function $f(\mathbf x,\bm\lambda)=e^{-\langle \mathbf x,\bm\lambda\rangle}$. Then we have
	\begin{eqnarray*}
		\ar\ar f_1'(\mathbf x,\bm\lambda)=-\lambda_1 f(\mathbf x,\bm\lambda),\quad f_2'(\mathbf x,\bm\lambda)=-\lambda_2f(\mathbf x,\bm\lambda),\quad f_3'(\mathbf x,\bm\lambda)=-x_1f(\mathbf x,\bm\lambda),\\
		\ar\ar f_4'(\mathbf x,\bm\lambda)=-x_2f(\mathbf x,\bm\lambda),\quad f_{11}''(\mathbf x,\bm\lambda)=\lambda_1^2f(\mathbf x,\bm\lambda),\quad\  f_{22}''(\mathbf x,\bm\lambda)=\lambda_2^2f(\mathbf x,\bm\lambda),\\
		\ar\ar f_{33}''(\mathbf x,\bm\lambda)=x_1^2f(\mathbf x,\bm\lambda),\quad \ f_{44}''(\mathbf x,\bm\lambda)=x_2^2f(\mathbf x,\bm\lambda),\quad \ f_{12}''(\mathbf x,\bm\lambda)=f_{21}''(\mathbf x,\bm\lambda)=\lambda_1\lambda_2f(\mathbf x,\bm\lambda).
	\end{eqnarray*}
	By It\^o's formula, for $t\geq r\geq 0$ and $\bm\lambda\in\mathbb R_+^2$,
	{\small
	\begin{eqnarray*}
		e^{-\langle\mathbf X(t),\bm\lambda\rangle}\ar=\ar e^{-\langle\mathbf X(r),\mathbf v_{r,t}(\bm\lambda)\rangle}+\int_r^t f_1'(\mathbf X(s-),\mathbf v_{s-,t}(\bm\lambda))\d X_1(s)\\
		\ar\ar+\int_r^t f_{12}''(\mathbf X(s-),\mathbf v_{s-,t}(\bm\lambda))\d\langle M_1,M_2\rangle(s)+\int_r^t f_2'(\mathbf X(s-),\mathbf v_{s-,t}(\bm\lambda))\d X_2(s)\\
		\ar\ar+\int_r^t f_3'(\mathbf X(s-),\mathbf v_{s-,t}(\bm\lambda))\d v_{1,s,t}(\bm\lambda)+\int_r^t f_4'(\mathbf X(s-),\mathbf v_{s-,t}(\bm\lambda))\d v_{2,s,t}(\bm\lambda)\\
		\ar\ar+\frac{1}{2}\int_r^t f_{11}''(\mathbf X(s-),\mathbf v_{s-,t}(\bm\lambda))\d C_1(s)+\frac{1}{2}\int_r^t f_{22}''(\mathbf X(s-),\mathbf v_{s-,t}(\bm\lambda))\d C_2(s)\\
		\ar\ar+\sum_{s\in(r,t]\backslash D}[f(\mathbf X(s),\mathbf v_{s,t}(\bm\lambda))-f(\mathbf X(s-),\mathbf v_{s-,t}(\bm\lambda))\\
		\ar\ar\quad\qquad\qquad-f_1'(\mathbf X(s-),\mathbf v_{s-,t}(\bm\lambda))\triangle X_1(s)-f_2'(\mathbf X(s-),\mathbf v_{s-,t}(\bm\lambda))\triangle X_2(s)]\\
		\ar\ar+\sum_{s\in(r,t]\cap D}[f(\mathbf X(s),\mathbf v_{s,t}(\bm\lambda))-f(\mathbf X(s-),\mathbf v_{s-,t}(\bm\lambda))\\
		\ar\ar\quad\qquad\qquad-f_1'(\mathbf X(s-),\mathbf v_{s-,t}(\bm\lambda))\triangle X_1(s)-f_3'(\mathbf X(s-),\mathbf v_{s-,t}(\bm\lambda))\triangle v_{1,s,t}(\bm\lambda)\\
		\ar\ar\quad\qquad\qquad-f_2'(\mathbf X(s-),\mathbf v_{s-,t}(\bm\lambda))\triangle X_2(s)-f_4'(\mathbf X(s-),\mathbf v_{s-,t}(\bm\lambda))\triangle v_{2,s,t}(\bm\lambda)]\\
		\ar=\ar e^{-\langle\mathbf X(r),\mathbf v_{r,t}(\bm\lambda)\rangle}+Z(t)-\sum_{i,j=1,2,~i\neq j}\int_0^r e^{-\langle\mathbf X(s-),\mathbf v_{s-,t}(\bm\lambda)\rangle}X_i(s-)v_{j,s,t}(\bm\lambda)\overline b_{ij}(\d s)\\
		\ar\ar+\sum_{i=1}^2\int_0^re^{-\langle\mathbf X(s-),\mathbf v_{s-,t}(\bm\lambda)\rangle}X_i(s-)v_{i,s,t}(\bm\lambda)b_{ii}(\d s) \\
		\ar\ar+\sum_{i=1}^2\int_0^r e^{-\langle\mathbf X(s-),\mathbf v_{s-,t}(\bm\lambda)\rangle}v_{i,s-,t}(\bm\lambda)\big[\d X_i^{mart.}(s)+A_i(\d s)]\\
		\ar\ar+\sum_{i=1}^2\int_0^r e^{-\langle\mathbf X(s-),\mathbf v_{s-,t}(\bm\lambda)\rangle}v_{i,s-,t}(\bm\lambda)^2\left[X_i(s-)c_i(\d s)-\frac{1}{2}C_i(\d s)\right]\\
		\ar\ar+\sum_{i=1}^2\int_0^r\int_{\mathbb R_+^2\backslash\{\mathbf 0\}}e^{-\langle\mathbf X(s-),\mathbf v_{s-,t}(\bm\lambda)\rangle}X_i(s-)K(\mathbf v_{s,t}(\bm\lambda),\mathbf z)m_i(\d s,\d \mathbf z)\\
		\ar\ar-\int_0^r e^{-\langle\mathbf X(s-),\mathbf v_{s-,t}(\bm\lambda)\rangle}v_{1,s-,t}(\bm\lambda)v_{2,s-,t}(\bm\lambda)\d\langle M_1,M_2\rangle(s)\\
		\ar\ar-\sum_{s\in(0,r]\backslash D}e^{-\langle\mathbf X(s-),\mathbf v_{s-,t}(\bm\lambda)\rangle}K(\mathbf v_{s,t}(\bm\lambda),\triangle\mathbf X(s))\\
		\ar\ar-\sum_{s\in(0,r]\cap D}\bigg[e^{-\langle\mathbf X(s),\mathbf v_{s,t}(\bm\lambda)\rangle}-e^{-\langle\mathbf X(s-),\mathbf v_{s-,t}(\bm\lambda)\rangle}\\
		\ar\ar\quad\qquad\qquad+e^{-\langle\mathbf X(s-),\mathbf v_{s-,t}(\bm\lambda)\rangle}\langle\mathbf v_{s-,t}(\bm\lambda),\triangle\mathbf X^{mart.}(s)\rangle\\
		\ar\ar\quad\qquad\qquad+e^{-\langle\mathbf X(s-),\mathbf v_{s-,t}(\bm\lambda)\rangle}v_{1,s-,t}(\bm\lambda)\triangle A_1(s)+e^{-\langle\mathbf X(s-),\mathbf v_{s-,t}(\bm\lambda)\rangle}v_{2,s-,t}(\bm\lambda)\triangle A_2(s)\\
		\ar\ar\quad\qquad\qquad+e^{-\langle\mathbf X(s-),\mathbf v_{s-,t}(\bm\lambda)\rangle}X_1(s-)\big(v_{1,s,t}(\bm\lambda)\triangle b_{11}(s)-v_{2,s,t}(\bm\lambda)\triangle\bar b_{12}(s)\big)\\
		\ar\ar\quad\qquad\qquad+e^{-\langle\mathbf X(s-),\mathbf v_{s-,t}(\bm\lambda)\rangle}X_2(s-)\big(v_{2,s,t}(\bm\lambda)\triangle b_{22}(s)-v_{1,s,t}(\bm\lambda)\triangle\bar b_{21}(s)\big)\\
		\ar\ar\quad\qquad\qquad+\int_{\mathbb R_+^2\backslash\{\mathbf 0\}}e^{-\langle\mathbf X(s-),\mathbf v_{s-,t}(\bm\lambda)\rangle}X_1(s-)K(\mathbf v_{s,t}(\bm\lambda),\mathbf z)m_1(\{s\},\d\mathbf z) \\
		\ar\ar\quad\qquad\qquad+\int_{\mathbb R_+^2\backslash\{\mathbf 0\}}e^{-\langle\mathbf X(s-),\mathbf v_{s-,t}(\bm\lambda)\rangle}X_2(s-)K(\mathbf v_{s,t}(\bm\lambda),\mathbf z)m_2(\{s\},\d\mathbf z)\bigg],
	\end{eqnarray*}
	}where $X_i^{mart.}(t):=X_i(t)-A_i(t)$ is a local martingale and $-Z(r)$ is defined by the right hand side of the second equality getting rid of $e^{-\langle\mathbf X(r),\mathbf v(r,t,\bm\lambda)\rangle}+Z(t)$. The second equality holds because of the continuity of $t\mapsto c_i(t)$ and $\mathbf v_{s,t}(\bm\lambda)=\mathbf v_{s-,t}(\bm\lambda)$ for $s\in (0,t]\backslash D$.
	Taking the conditional expectation in above equation, we have
	\begin{eqnarray*}
		\mathbf E[-Z(t)|\mathscr F_r]\ar=\ar\sum_{i,j=1,2,~i\neq j}\!\int_{(0,r]\backslash D}\!\!e^{-\langle\mathbf X(s-),\mathbf v_{s-,t}(\bm\lambda)\rangle}v_{i,s-,t}(\bm\lambda)[A_i(\d s)\!-\!X_{j}(s-)\bar b_{ji}(\d s)\!+\!X_i(s-)b_{ii}(\d s)]\\
		\ar\ar+\sum_{i=1}^2\int_0^r e^{-\langle\mathbf X(s-),\mathbf v_{s-,t}(\bm\lambda)\rangle}v_{i,s-,t}(\bm\lambda)^2\left[X_i(s-)c_i(\d s)-\frac{1}{2}C_i(\d s)\right]\\
		\ar\ar-\int_0^r e^{-\langle\mathbf X(s-),\mathbf v_{s-,t}(\bm\lambda)\rangle}v_{1,s-,t}(\bm\lambda)v_{2,s-,t}(\bm\lambda)\d\langle M_1,M_2\rangle(s)+\text{mart.}\\
		\ar\ar+\sum_{i=1}^2\int_{(0,r]\backslash D}\int_{\mathbb R_+^2\backslash\{\mathbf 0\}}e^{-\langle\mathbf X(s-),\mathbf v_{s-,t}(\bm\lambda)\rangle}X_i(s-)K(\mathbf v_{s,t}(\bm\lambda),\mathbf z)m_i(\d s,\d \mathbf z)\\
		\ar\ar-\sum_{s\in(0,r]\backslash D}e^{-\langle\mathbf X(s-),\mathbf v_{s-,t}(\bm\lambda)\rangle}K(\mathbf v_{s,t}(\bm\lambda),\triangle\mathbf X(s))\\
		\ar\ar-\sum_{s\in(0,r]\cap D}e^{-\langle\mathbf X(s-),\mathbf v_{s-,t}(\bm\lambda)\rangle}\langle\mathbf v_{s-,t}(\bm\lambda),\triangle\mathbf X^{mart.}(s)\rangle.
	\end{eqnarray*}
	The first equality holds because of the absolute continuity property of $b_{ii,c},\bar b_{ij,c},i=1,2,j\neq i$ with respect to Lebesgue measure. Then by the uniqueness of canonical decompositions of semimartingales, we have
	\begin{eqnarray*}
		\ar\ar \mathbf 1_{D^c}(s)A_i(\d s)=\mathbf 1_{D^c}(s)[X_j(s-)\bar b_{ji}(\d s)-X_i(s-)b_{ii}(\d s)],\\
		\ar\ar C_i(\d s)=2X_i(s-)c_i(\d s)=2c_i(\d s)\int_0^\infty \mathbf 1_{\{u\leq X_i(s-)\}}\d u,\\
		\ar\ar\d \langle M_1,M_2\rangle(s)=0
	\end{eqnarray*}
	and
	\begin{eqnarray*}
		\mathbf 1_{D^c}(s)\hat N_0(\d s,\d \mathbf z)\ar=\ar\mathbf 1_{D^c}(s)[X_1(s-)m_1(\d s,\d \mathbf z)+X_2(s-)m_2(\d s,\d \mathbf z)]\\
			\ar=\ar\mathbf 1_{D^c}(s)\int_0^\infty [\mathbf 1_{\{u\leq X_1(s-)\}} m_1(\d s,\d \mathbf z)+\mathbf 1_{\{u\leq X_2(s-)\}} m_2(\d s,\d \mathbf z)]\d u,
	\end{eqnarray*}
	thus $M_1$ and $M_2$ are independent; see, e.g., Karatzas and Shreve (1988), Theorem 4.13. By Theorem III.6 in El Karoui and M\'el\'eard (1990), on the extension of the original probability space there exists two independent Gaussian white noises $W_1(\d s,\d u)$ and $W_2(\d s,\d u)$ on $(0,\infty)^2$ with intensity $2c_1(\d s)\d u$ and $2c_2(\d s)\d u$, respectively, such that
	\begin{equation*}
		M_i(t)=\int_0^t\int_0^\infty \mathbf 1_{\{u\leq X_i(s-)\}}W_i(\d s,\d u)=\int_0^t\int_0^{X_i(s-)}W_i(\d s,\d u).
	\end{equation*}
	Moreover, by Ikeda and Watanabe (1989), on a further extension of the original probability space we can define two independent Poisson random measures $M_{1,D^c}(\d s,\d \mathbf z,\d u)$ and $M_{2,D^c}(\d s,\d \mathbf z,$ $\d u)$ with intensity $\mathbf 1_{D^c}(s)m_1(\d s,\d \mathbf z)\d u$ and $\mathbf 1_{D^c}(s)m_2(\d s,\d \mathbf z)\d u$ such that
	\begin{eqnarray*}
		\int_0^t\int_{\mathbb R_+^2\backslash\{\mathbf 0\}}z_i\mathbf 1_{D^c}(s)\tilde N_0(\d s,\d \mathbf z)\ar=\ar\int_0^t\int_{\mathbb R_+^2\backslash\{\mathbf 0\}}\int_0^{X_1(s-)}z_i\tilde M_{1,D^c}(\d s,\d \mathbf z,\d u)\\
		\ar\ar+\int_0^t\int_{\mathbb R_+^2\backslash\{\mathbf 0\}}\int_0^{X_2(s-)}z_i\tilde M_{2,D^c}(\d s,\d \mathbf z,\d u).
	\end{eqnarray*}
	Furthermore, from \eqref{X_i(t)decomposition},  we see that, for $s\in D$, the process a.s. makes a jump as
	\begin{equation*}
		\triangle X_i(s)=\triangle A_i(s)+\int_{\mathbb R^2}z_i\tilde N_0(\{s\},\d \mathbf z).
	\end{equation*}
	By Proposition \ref{TCBVE realization}, we have
	\begin{equation*}
		\mathbf E(e^{-\langle\bm\lambda,\triangle\mathbf X(s)\rangle}|\mathscr F_{s-})=e^{(\lambda_1-v_{1,s-,s}(\bm\lambda))X_1(s-)+(\lambda_2-v_{2,s-,s}(\bm\lambda))X_2(s-)},\quad \bm\lambda\in\mathbb R_+^2,
	\end{equation*}
	where
	\begin{equation*}
		\lambda_i-v_{i,s-,s}(\bm\lambda)=\lambda_i\triangle b_{ii}(s)-\lambda_j\triangle\overline b_{ij}(s)+\int_{\mathbb R_+^2\backslash\{\mathbf 0\}}K(\bm\lambda,\mathbf z)m_i(\{s\},\d \mathbf z).
	\end{equation*}
	Then by Proposition \ref{composite Levy-Ito representation}, we can make another extension of the probability space and define two independent Poisson random measures $N_{1,s}(\d \mathbf z,\d u)$ and $N_{2,s}(\d \mathbf z,\d u)$ with intensity $m_1(\{s\},\d \mathbf z)\d u$ and $m_2(\{s\},\d \mathbf z)\d u$ such that $N_{1,s}$ and $N_{2,s}$ are independent of $\mathscr F_{s-}$ and
	{\small
	\begin{equation*}
		\triangle X_i(s)=X_j(s-)\triangle\bar b_{ji}(s)-X_i(s-)\triangle b_{ii}(s)+\int_{\mathbb R_+^2\backslash\{\mathbf 0\}}\int_0^{X_1(s-)}z_i\tilde N_{1,s}(\d \mathbf z,\d u)+\int_{\mathbb R_+^2\backslash\{\mathbf 0\}}\int_0^{X_2(s-)}z_i\tilde N_{2,s}(\d \mathbf z,\d u).
	\end{equation*}
	}We see that $\triangle A_i(s)=X_j(s-)\triangle\bar b_{ji}(s)-X_i(s-)\triangle b_{ii}(s)$ for $s\in D$, thus $A_i(\d s)=X_j(s-)\bar b_{ji}(\d s)-X_i(s-)b_{ii}(\d s)$ for every $s\geq 0$. Let $M_{1,D}(\d s,\d \mathbf z,\d u)$ and $M_{2,D}(\d s,\d \mathbf z,\d u)$ be independent random measures on $(0,\infty)\times (\mathbb R_+^2\backslash\{\mathbf 0\})\times(0,\infty)$ defined by
	\begin{equation*}
		M_{i,D}((0,t]\times B)=\sum_{s\in (0,t]\cap D} N_{i,s}(B),\quad t\geq 0,\ B\in\mathscr B((\mathbb R_+^2\backslash\{\mathbf 0\})\times(0,\infty)).
	\end{equation*}
	Then $M_{i,D}(\d s,\d \mathbf z,\d u)$ is a Poisson random measure with intensity $\mathbf 1_D(s)m_i(\d s,\d\mathbf z)\d u$. By the decomposition representation \eqref{X_i(t)decomposition}, we see that $\{\mathbf X(t):t\geq 0\}$ is a solution to \eqref{4 SIE}. The independence of the white noises and  Poisson noises follows by standard arguments; see, e.g., Ikeda and Watanabe (1989), pages 77-78. \eproof

\section{General stochastic equations}

\setcounter{equation}{0}

In this section, we prove that the TCBVE-process defined by \eqref{Laplace transform} and \eqref{backward eq} is the pathwise unique solution to \eqref{SIE}.

\medskip\noindent\textit{Proof of Theorem~\ref{SIE solution}.~} \textit{Step 1.} We first consider the case under condition \eqref{additional integrability condition}. By Proposition \ref{SIE's weak solution} and \ref{comparison property}, we see that \eqref{4 SIE} has a pathwise unique weak solution, which is a TCBVE-process with transition semigroup $(Q_{r,t})_{t\geq r}$ defined by \eqref{Laplace transform} and \eqref{backward eq}. Then by standard arguments, we see that \eqref{4 SIE} has a pathwise unique solution; see, e.g., Situ (2005), page 104.
	
\textit{Step 2.} Now we consider the general case. For simplicity, let $(i,j)=(1,2)$ or $(2,1)$. For $k\geq 1$, there is a pathwise unique solution $\{\mathbf X^{(k)}(t):t\geq 0\}$ to the stochastic equations
	\begin{equation}\small\label{modified X_i(t)}
	\begin{split}
		&\!X_i(t)=X_i(0)-\int_0^t X_i(s-)b_{ii}(\d s)+\int_0^t X_j(s-)\bar b_{ji}(\d s)+\int_0^t\int_0^{X_i(s-)}W_i(\d s,\d u)\\
		&\!+\!\int_0^t\int_{\mathbb R_+^2\backslash\{\mathbf 0\}}\int_0^{X_i(s-)}(z_i\wedge k)\tilde M_i(\d s,\d \mathbf z,\d u)-\int_0^t\int_{\{z_i>k,z_j\geq 0\}}(z_i-k)X_i(s-)m_i(\d s,\d \mathbf z)\\
		&\!+\!\int_0^t\int_{\mathbb R_+^2\backslash\{\mathbf 0\}}\int_0^{X_j(s-)}(z_i\wedge k)\tilde M_j(\d s,\d \mathbf z,\d u)-\int_0^t\int_{\{z_i>k,z_j\geq 0\}}(z_i-k)X_j(s-)m_j(\d s,\d \mathbf z)	
	\end{split}
	\end{equation}
	and the solution is a TCBVE-process with cumulant semigroup $(\mathbf v^{(k)}_{r,t}(\bm\lambda))_{t\geq r}$ defined by
	\begin{equation}\small\label{modified cumulant semigroup}
		\begin{split}
			&v_{i,r,t}(\bm\lambda)=\lambda_i-\int_r^t v_{i,s,t}(\bm\lambda)b_{ii}(\d s)+\int_r^t v_{j,s,t}(\bm\lambda)\bar b_{ij}(\d s)-\int_r^t v_{i,s,t}(\bm\lambda)^2c_i(\d s)\\
			&\quad-\int_r^t\int_{\{z_i>k,z_j\geq 0\}}v_{i,s,t}(\bm\lambda)(z_i-k)m_i(\d s,\d \mathbf z)-\int_r^t\int_{\{z_j>k,z_i\geq 0\}}v_{j,s,t}(\bm\lambda)(z_j-k)m_i(\d s,\d \mathbf z)\\
			&\quad-\int_r^t\int_{\mathbb R_+^2\backslash\{\mathbf 0\}}K(\mathbf v_{s,t}(\bm\lambda),z_1\wedge k, z_2\wedge k)m_i(\d s,\d \mathbf z).
		\end{split}
	\end{equation}
	We can rewrite the above stochastic equation into the equivalent form:
	\begin{equation}\small\label{modified X_i(t) equivalent form}
	\begin{split}
		&X_i(t)=X_i(0)-\int_0^t X_i(s-)b_{ii}(\d s)+\int_0^t X_j(s-)\bar b_{ji}(\d s)+\int_0^t\int_0^{X_i(s-)}W_i(\d s,\d u)\\
		&+\int_0^t\int_{\{0\leq z_i\leq k,z_j\geq 0\}\backslash\{\mathbf 0\}}\int_0^{X_i(s-)}z_i\tilde M_i(\d s,\d \mathbf z,\d u)-\int_0^t\int_{\{z_i>k,z_j\geq 0\}}z_iX_i(s-)m_i(\d s,\d \mathbf z)\\
		&+\int_0^t\int_{\{0\leq z_i\leq k,z_j\geq 0\}\backslash\{\mathbf 0\}}\int_0^{X_j(s-)}z_i\tilde M_j(\d s,\d \mathbf z,\d u)-\int_0^t\int_{\{z_i>k,z_j\geq 0\}}z_iX_j(s-)m_j(\d s,\d \mathbf z)\\
		&+\int_0^t\int_{\{z_i>k,z_j\geq 0\}}\int_0^{X_i(s-)}k M_i(\d s,\d \mathbf z,\d u)+\int_0^t\int_{\{z_i>k,z_j\geq 0\}}\int_0^{X_j(s-)}k M_j(\d s,\d \mathbf z,\d u).
		\end{split}
	\end{equation}
	This is a special case of \eqref{4 SIE} with $z_i\mathbf 1_{\{z_i>k,z_j\geq 0\}}(\mbf z)M_p(\d s,\d \mathbf z,\d u)$ replaced by $k\mathbf 1_{\{z_i>k,z_j\geq 0\}}(\mbf z)$\ $M_p(\d s,\d \mathbf z,\d u)$ for $p=1,2$, which modifies the magnitudes of the large jumps. Let $\zeta_{0,k}=0$ and for $l\geq 0$, define
	\begin{eqnarray*}
		\zeta_{l+1,k}=\inf\bigg\{t>0:\int_{\zeta_{l,k}}^{\zeta_{l,k}+t}\int_{\mathbb R_+^2\backslash [0,k]^2}\int_0^{X_1^{(k)}(s-)}M_1(\d s,\d \mathbf z,\d u)\\
		+\int_{\zeta_{l,k}}^{\zeta_{l,k}+t}\int_{\mathbb R_+^2\backslash [0,k]^2}\int_0^{X_2^{(k)}(s-)}M_2(\d s,\d \mathbf z,\d u)\geq 1\bigg\},
	\end{eqnarray*}
	inductively. Then $\eta_{n,k}:=\sum_{l=1}^n\zeta_{l,k}$ is the time when $\{\mathbf X^{(k)}:t\geq 0\}$ has the $n$th modified jump. By the c\`adl\`ag property of $\{\mathbf X^{(k)}:t\geq 0\}$, we have $\eta_{n,k}\to\infty$ as $n\to\infty$. Clearly, we see that the modified solution is same as the original one in the time interval $[0,\eta_{1,k})$, and then $\{X^{(k)}(t):t\geq 0\}$ takes its first modified jump at time $\eta_{1,k}=\zeta_{1,k}$. For $k\geq 1$, we have $\eta_{1,k}\leq \eta_{1,k+1}$, which implies $X_i^{(k)}(t)=X_i^{(k+1)}(t)$ for $0\leq t<\eta_{1,k}$. By \eqref{modified X_i(t) equivalent form}, we have
	{\small
	\begin{eqnarray*}
		\ar\ar\!\!\!\!\triangle X_i^{(k)}(\eta_{1,k})= X_j^{(k)}(\eta_{1,k}-)\triangle\bar b_{ji}(\eta_{1,k})-X_i^{(k)}(\eta_{1,k}-)\triangle b_{ii}(\eta_{1,k})\\
		\ar\ar\!\!\!+\int_{\{0\leq z_i\leq k,z_j\geq 0\}\backslash\{\mathbf 0\}}\int_0^{X_i^{(k)}(\eta_{1,k}-)}z_i\tilde M_i(\{\eta_{1,k}\},\d \mathbf z,\d u)\\
		\ar\ar\!\!\!+\int_{\{0\leq z_i\leq k,z_j\geq 0\}\backslash\{\mathbf 0\}}\int_0^{X_j^{(k)}(\eta_{1,k}-)}z_i\tilde M_j(\{\eta_{1,k}\},\d \mathbf z,\d u)\\
		\ar\ar\!\!\!+\int_{\{z_i>k,z_j\geq 0\}}\int_0^{X_i^{(k)}(\eta_{1,k}-)}k M_i(\{\eta_{1,k}\},\d \mathbf z,\d u)+\int_{\{z_i>k,z_j\geq 0\}}\int_0^{X_j^{(k)}(\eta_{1,k}-)}k M_j(\{\eta_{1,k}\},\d \mathbf z,\d u)\\
		\ar\ar\!\!\!-\int_{\{z_i>k,z_j\geq 0\}}z_iX_i^{(k)}(\eta_{1,k}-)m_i(\{\eta_{1,k}\},\d \mathbf z)-\int_{\{z_i>k,z_j\geq 0\}}z_iX_j^{(k)}(\eta_{1,k}-)m_j(\{\eta_{1,k}\},\d \mathbf z)
	\end{eqnarray*}
	}and
	{\small
	\begin{eqnarray*}
		\triangle X_i^{(k+1)}(\eta_{1,k})\ar=\ar X_j^{(k)}(\eta_{1,k}-)\triangle\bar b_{ji}(\eta_{1,k})-X_i^{(k)}(\eta_{1,k}-)\triangle b_{ii}(\eta_{1,k})\\
		\ar\ar+\int_{\{0\leq z_i\leq k+1,z_j\geq 0\}\backslash\{\mathbf 0\}}\int_0^{X_i^{(k)}(\eta_{1,k}-)}z_i\tilde M_i(\{\eta_{1,k}\},\d \mathbf z,\d u)\\
		\ar\ar+\int_{\{0\leq z_i\leq k+1,z_j\geq 0\}\backslash\{\mathbf 0\}}\int_0^{X_j^{(k)}(\eta_{1,k}-)}z_i\tilde M_j(\{\eta_{1,k}\},\d \mathbf z,\d u)\\
		\ar\ar+\int_{\{z_i>k+1,z_j\geq 0\}}\int_0^{X_i^{(k)}(\eta_{1,k}-)}(k+1) M_i(\{\eta_{1,k}\},\d \mathbf z,\d u)\\
		\ar\ar+\int_{\{z_i>k+1,z_j\geq 0\}}\int_0^{X_j^{(k)}(\eta_{1,k}-)}(k+1) M_j(\{\eta_{1,k}\},\d \mathbf z,\d u)\\
		\ar\ar-\int_{\{z_i>k+1,z_j\geq 0\}}z_iX_i^{(k)}(\eta_{1,k}-)m_i(\{\eta_{1,k}\},\d \mathbf z)\\
		\ar\ar-\int_{\{z_i>k+1,z_j\geq 0\}}z_iX_j^{(k)}(\eta_{1,k}-)m_j(\{\eta_{1,k}\},\d \mathbf z)\\
		\ar=\ar X_j^{(k)}(\eta_{1,k}-)\triangle\bar b_{ji}(\eta_{1,k})-X_i^{(k)}(\eta_{1,k}-)\triangle b_{ii}(\eta_{1,k})\\
		\ar\ar+\int_{\{0\leq z_i\leq k,z_j\geq 0\}\backslash\{\mathbf 0\}}\int_0^{X_i^{(k)}(\eta_{1,k}-)}z_i\tilde M_i(\{\eta_{1,k}\},\d \mathbf z,\d u)\\
		\ar\ar+\int_{\{0\leq z_i\leq k,z_j\geq 0\}\backslash\{\mathbf 0\}}\int_0^{X_j^{(k)}(\eta_{1,k}-)}z_i\tilde M_j(\{\eta_{1,k}\},\d \mathbf z,\d u)\\
		\ar\ar+\int_{\{z_i>k+1,z_j\geq 0\}}\int_0^{X_i^{(k)}(\eta_{1,k}-)}(k+1) M_i(\{\eta_{1,k}\},\d \mathbf z,\d u)\\
		\ar\ar+\int_{\{z_i>k+1,z_j\geq 0\}}\int_0^{X_j^{(k)}(\eta_{1,k}-)}(k+1) M_j(\{\eta_{1,k}\},\d \mathbf z,\d u)\\
		\ar\ar-\int_{\{z_i>k,z_j\geq 0\}}z_iX_i^{(k)}(\eta_{1,k}-)m_i(\{\eta_{1,k}\},\d \mathbf z)\\
		\ar\ar-\int_{\{z_i>k,z_j\geq 0\}}z_iX_j^{(k)}(\eta_{1,k}-)m_j(\{\eta_{1,k}\},\d \mathbf z)\\
		\ar\ar+\int_{\{k< z_i\leq k+1,z_j\geq 0\}\backslash\{\mathbf 0\}}\int_0^{X_i^{(k)}(\eta_{1,k}-)}z_i M_i(\{\eta_{1,k}\},\d \mathbf z,\d u)\\
		\ar\ar+\int_{\{k<z_i\leq k+1,z_j\geq 0\}\backslash\{\mathbf 0\}}\int_0^{X_j^{(k)}(\eta_{1,k}-)}z_i M_j(\{\eta_{1,k}\},\d \mathbf z,\d u).
	\end{eqnarray*}
	}Then $\triangle X_i^{(k)}(\eta_{1,k})\leq \triangle X_i^{(k+1)}(\eta_{1,k})$, thus $X_i^{(k)}(\eta_{1,k})\leq X_i^{(k+1)}(\eta_{1,k})$. By applying Proposition \ref{comparison property} successively at the stopping times $\eta_{n,k},n\geq 1$, we have $X_i^{(k)}(t)\leq X_i^{(k+1)}(t)$ holds for all $t\geq 0$ and $k\geq 1$. For $r\geq 0,\mathbf x\in\mathbb R_+^2$, let $\{\mathbf X^{(k)}(\mathbf x,r,t):t\geq r\}$ be the pathwise unique solution to stochastic equations
	{\small
	\begin{eqnarray*}
		X_i(t)\ar=\ar x_i-\int_r^t X_i(s-)b_{ii}(\d s)+\int_r^t X_j(s-)\bar b_{ji}(\d s)+\int_r^t\int_0^{X_i(s-)}W_i(\d s,\d u)\\
		\ar\ar+\int_r^t\int_{\mathbb R_+^2\backslash\{\mathbf 0\}}\int_0^{X_i(s-)}(z_i\wedge k)\tilde M_i(\d s,\d \mathbf z,\d u)-\int_r^t\int_{\{z_i>k,z_j\geq 0\}}(z_i-k)X_i(s-)m_i(\d s,\d \mathbf z)\\
		\ar\ar+\int_r^t\int_{\mathbb R_+^2\backslash\{\mathbf 0\}}\int_0^{X_j(s-)}(z_i\wedge k)\tilde M_j(\d s,\d \mathbf z,\d u)-\int_r^t\int_{\{z_i>k,z_j\geq 0\}}(z_i-k)X_j(s-)m_j(\d s,\d \mathbf z).
	\end{eqnarray*}
	}Then we have $X_i^{(k+1)}(\mathbf x,r,t)\geq X_i^{(k)}(\mathbf x,r,t)$ and
	\begin{eqnarray*}
		v_{i,r,t}^{(k+1)}(\bm\lambda)\ar=\ar -\log\mathbf E\exp\{-\langle\bm\lambda,\mathbf X^{(k+1)}(\mathbf e_i,r,t)\rangle\}\geq-\log\mathbf E\exp\{-\langle\bm\lambda,\mathbf X^{(k)}(\mathbf e_i,r,t)\rangle\}=v_{i,r,t}^{(k)}(\bm\lambda).
	\end{eqnarray*}
	By upper estimation in Li and Zhang (2025+), we have $v_{i,r,t}^{(k)}(\bm\lambda)\leq U_i(0,t,\bm\lambda)$. Moreover, let $k\to\infty$ in \eqref{modified cumulant semigroup}, we see that the limit $r\mapsto \mathbf v_{r,t}(\bm\lambda):=\uparrow\lim_{k\to\infty} \mathbf v^{(k)}_{r,t}(\bm\lambda)$ is the unique bounded positive solution to \eqref{backward eq}. By Proposition \ref{EX_(t) estimate}, we have $\lim_{k\to\infty}\eta_{1,k}=\lim_{k\to\infty}\zeta_{1,k}\geq\lim_{k\to\infty}\tau_{k,1}\wedge\tau_{k,2}=\infty$. Then we can construct a c\`adl\`ag process $\{\mathbf Y(t):t\geq 0\}$ such that $\mathbf Y(t)=\mathbf X^{(k)}(t)$ for $0\leq t< \eta_{1,k}$ and $k\geq 1$. Furthermore, for $0\leq t<\eta_{1,k}$ we have
	\begin{equation*}\small
		\begin{split}
			&\int_0^t\int_{\{z_i>k,z_j\geq 0\}}\int_0^{X_i^{(k)}(s-)}z_i M_i(\d s,\d \mathbf z,\d u)=\int_0^t\int_{\{z_i>k,z_j\geq 0\}}\int_0^{X_j^{(k)}(s-)}z_i M_j(\d s,\d \mathbf z,\d u)\\
			&=\int_0^t\int_{\{z_i>k,z_j\geq 0\}}\int_0^{X_i^{(k)}(s-)}k M_i(\d s,\d \mathbf z,\d u)=\int_0^t\int_{\{z_i>k,z_j\geq 0\}}\int_0^{X_j^{(k)}(s-)}k M_j(\d s,\d \mathbf z,\d u)=0,
		\end{split}
	\end{equation*}
	so \eqref{modified X_i(t) equivalent form} implies
	{\small
	\begin{eqnarray*}
		Y_i(t)\ar=\ar X_i(0)-\int_0^t Y_i(s-)b_{ii}(\d s)+\int_0^t Y_j(s-)\bar b_{ji}(\d s)+\int_0^t\int_0^{Y_i(s-)}W_i(\d s,\d u)\\
		\ar\ar+\int_0^t\int_{\{0\leq z_i\leq k,z_j\geq 0\}\backslash\{\mathbf 0\}}\int_0^{Y_i(s-)}z_i\tilde M_i(\d s,\d \mathbf z,\d u)\\
		\ar\ar+\int_0^t\int_{\{0\leq z_i\leq k,z_j\geq 0\}\backslash\{\mathbf 0\}}\int_0^{Y_j(s-)}z_i\tilde M_j(\d s,\d \mathbf z,\d u)\\
		\ar\ar+\int_0^t\int_{\{z_i>k,z_j\geq 0\}}\int_0^{Y_i(s-)}z_i M_i(\d s,\d \mathbf z,\d u)+\int_0^t\int_{\{z_i>k,z_j\geq 0\}}\int_0^{Y_j(s-)}z_i M_j(\d s,\d \mathbf z,\d u)\\
		\ar\ar-\int_0^t\int_{\{z_i>k,z_j\geq 0\}}z_iY_i(s-)m_i(\d s,\d \mathbf z)-\int_0^t\int_{\{z_i>k,z_j\geq 0\}}z_iY_j(s-)m_j(\d s,\d \mathbf z)\\
		\ar=\ar X_i(0)-\int_0^t Y_i(s-)b_{ii}(\d s)+\int_0^t Y_j(s-)\bar b_{ji}(\d s)+\int_0^t\int_0^{Y_i(s-)}W_i(\d s,\d u)\\
		\ar\ar+\int_0^t\int_{\mathbb R_+^2\backslash\{\mathbf 0\}}\int_0^{Y_i(s-)}z_i\tilde M_i(\d s,\d \mathbf z,\d u)+\int_0^t\int_{\mathbb R_+^2\backslash\{\mathbf 0\}}\int_0^{Y_j(s-)}z_i\tilde M_j(\d s,\d \mathbf z,\d u).
 \end{eqnarray*}
}Since $\lim_{k\to\infty}\eta_{1,k}=\infty$, we can easily obtain $\{\mathbf Y(t):t\geq 0\}$ is a solution to \eqref{4 SIE}. On the other hand, it is clear that $\{X_i^{(k)}(t):t\geq 0\}$ increases to $\{Y_i(t):t\geq 0\}$ as $k\to\infty$, so $\{\mathbf Y(t):t\geq 0\}$ is a TCBVE-process with cumulant semigroup $(\mathbf v(r,t,\bm\lambda))_{t\geq r}$. Similarly, if $\{\mathbf Z(t):t\geq 0\}$ is also a solution to \eqref{4 SIE}, we have $\mathbf Z(t)=\mathbf X^{(k)}(t)$ for $0\leq t<\eta_{1,k}$ and $k\geq 1$, implying the pathwise uniqueness of solution to \eqref{4 SIE}.	\qed

Moreover, we can give a simple characterization of the distribution of extinction time.

\begin{corollary} Let $\tau_{\mathbf 0}=\inf\{t\geq 0: \mathbf X(t)=\mathbf 0 \}$. Then in the setting of Theorem \ref{SIE solution}, for any $t\geq 0$ we have
 \begin{equation}
\mathbf P\{\tau_{\mathbf 0}\leq t\}=\mathbf P\{\mathbf X(t)=\mathbf 0\}=\mathbf E(e^{-\langle \mathbf X(0),\mathbf v_{0,t}(\infty)\rangle}),
 \end{equation}
where $\mathbf v_{0,t}(\infty)=\lim_{\lambda_1,\lambda_2\uparrow\infty} \mathbf v_{0,t}(\bm\lambda)$. \end{corollary}

\section{Positive integral functionals}

In this section, we give characterizations of a class of positive integral functionals of the TCBVE-process in terms of Laplace transforms.

\begin{proposition}\label{positive integral functionals prop} Let $\mathbf X=(\Omega,\mathscr F,\mathscr F_{r,t},\mathbf X(t),\mathbf P_{r,\mathbf x})$ be a TCBVE-process started from time $r\ge 0$. For $i=1,2$, let $\zeta_i$ be a $\sigma$-finite measure on $\mathbb R_+$ satisfying $\zeta_i(B)<\infty$ for every bounded Borel set $B\subset \mathbb R_+$. Then for any $t\geq r$ and $\bm\lambda,\mbf x\in\mathbb R_+^2$ we have
 \begin{equation*}
\mathbf P_{r,\mathbf x}\exp\bigg\{-\langle\bm\lambda,\mathbf X(t)\rangle-\sum_{i=1}^2\int_r^t X_i(s)\zeta_i(\d s)\bigg\}=\exp\{-\langle\mathbf x,\mathbf u_{r,t}\rangle\},
 \end{equation*}
where $r\mapsto\mathbf u_{r,t}$ is the unique positive solution to
 \begin{equation}\label{u_r,t}
 \begin{split}
u_{i,r,t}=\lambda_i &+\int_r^t\zeta_i(\d s)-\int_r^t \big[u_{i,s,t}+\zeta_i(\{s\})\big]b_{ii}(\d s)+\int_r^t \big[u_{j,s,t}+\zeta_j(\{s\})\big]b_{ij}(\d s) \\
 &
-\int_r^t u_{i,s,t}^2 c_{i}(\d s)-\int_r^t\int_{\mathbb R_+^2\backslash\{\mathbf 0\}} K_i(\mathbf u_{s,t}+\bm\zeta(\{s\}),\mathbf z)m_i(\d s,\d \mathbf z),
 \end{split}
 \end{equation}
$r\in [0,t],~ i,j=1,2,~ j\neq i$ and $\bm\zeta(\{s\})=(\zeta_1(\{s\}),\zeta_2(\{s\}))$. \end{proposition}

\bproof	The uniqueness of solutions to \eqref{u_r,t} is obvious by Gronwall's inequality. Let $r=t$, we have $\mathbf u_{t,t}=\bm\lambda$. Next we will prove $\mathbf u_{r,t}$ is a solution to \eqref{u_r,t}. Define $\mathbf Z_r(s)$ by
 $$
Z_{i,r}(s)=\int_r^s X_i(h)\zeta_i(\d h),\quad s\geq r\geq 0.
 $$
Define a function $G(\mathbf u,\mathbf x,\mathbf z)=\exp\{-\langle\mathbf u,\mathbf x\rangle-z_1-z_2\}$ on $\mathbb R_+^6$. By simple calculus, we have
	\begin{eqnarray*}
	\begin{aligned}
		&G_1'(\mathbf u,\mathbf x,\mathbf z)=-x_1G(\mathbf u,\mathbf x,\mathbf z),\  G_2'(\mathbf u,\mathbf x,\mathbf z)=-x_2G(\mathbf u,\mathbf x,\mathbf z),\ G_3'(\mathbf u,\mathbf x,\mathbf z)=-u_1G(\mathbf u,\mathbf x,\mathbf z),\\
		&G_4'(\mathbf u,\mathbf x,\mathbf z)=-u_2G(\mathbf u,\mathbf x,\mathbf z),\ G_5'(\mathbf u,\mathbf x,\mathbf z)=G_6'(\mathbf u,\mathbf x,\mathbf z)=-G(\mathbf u,\mathbf x,\mathbf z),\\
		&G_{33}''(\mathbf u,\mathbf x,\mathbf z)=u_1^2G(\mathbf u,\mathbf x,\mathbf z),\ \ G_{44}''(\mathbf u,\mathbf x,\mathbf z)=u_2^2G(\mathbf u,\mathbf x,\mathbf z).
	\end{aligned}		
	\end{eqnarray*}
	Let $A:=J_{b_{11}}\cup J_{b_{22}}\cup J_{\bar b_{12}}\cup J_{\bar b_{21}}\cup J_{m_1}\cup J_{m_2}\cup J_{\zeta_1}\cup J_{\zeta_2}$. It is obvious that $A$ contains all discontinuous points of $\mathbf u_{r,t}$. For $s\in (r,t]\cap A$,
	{\small
	\begin{eqnarray*}
		\ar\ar\mathbf P_{r,\mathbf x}\exp\{-\langle\mathbf u_{s-,t},\mathbf X(s-)\rangle-Z_{1,r}(s-)-Z_{2,r}(s-)\}\\
		\ar\ar=\mathbf P_{r,\mathbf x}\exp\{-\langle\mathbf u_{s,t},\mathbf X(s)\rangle-Z_{1,r}(s)-Z_{2,r}(s)\}\\
		\ar\ar=\mathbf P_{r,\mathbf x}[\mathbf P_{r,\mathbf x}(\exp\{-\langle\mathbf u_{s,t},\mathbf X(s)\rangle-Z_{1,r}(s)-Z_{2,r}(s)\}|\mathscr F_{s-})]\\
		\ar\ar=\mathbf P_{r,\mathbf x}[\exp\{-\langle\mathbf u_{s,t}+\bm\zeta(\{s\}),\mathbf X(s-)\rangle-Z_{1,r}(s-)-Z_{2,r}(s-)\}\mathbf P_{r,\mathbf x}(\exp\{-\langle\mathbf u_{s,t}+\bm\zeta(\{s\}),\triangle\mathbf X(s-)\rangle\}|\mathscr F_{s-})]\\
		\ar\ar=\mathbf P_{r,\mathbf x}\Bigg[\exp\left\{-\langle\mathbf u_{s,t}+\bm\zeta(\{s\}),\mathbf X(s-)\rangle+\sum_{i=1}^2 \int_{\mathbb R_+^2\backslash\{\mathbf 0\}}X_i(s-)K(\mathbf u_{s,t}+\bm\zeta(\{s\}),\mathbf z)m_i(\{s\},\d \mathbf z)\right\}\\
		\ar\ar\qquad\cdot\exp\left\{-Z_{1,r}(s-)-Z_{2,r}(s-)+\sum_{i=1}^2 X_i(s-)\Big[(u_{i,s,t}+\zeta_i(\{s\}))\triangle b_{ii}(s)-(u_{j,s,t}+\zeta_j(\{s\}))\triangle\bar b_{ij}(s)\Big]\right\}\Bigg],
	\end{eqnarray*}
	}where we use Proposition \ref{TCBVE realization} for the last equality. Then we have
	\begin{eqnarray*}
		\triangle u_{i,s,t}=-\zeta_i(\{s\})\ar+\ar(u_{i,s,t}+\zeta_i(\{s\}))\triangle b_{ii}(s)-(u_{j,s,t}+\zeta_j(\{s\}))\triangle\bar b_{ij}(s)\\
		\ar+\ar\int_{\mathbb R_+^2\backslash\{\mathbf 0\}}K(\mathbf u_{s,t}+\bm\zeta(\{s\}),\mathbf z)m_i(\{s\},\d \mathbf z).
	\end{eqnarray*}
	By the Markov property of $X(t)$, we can easily obtain $\{e^{-\langle\mathbf u_{s,t},\mathbf X(s)\rangle-Z_{1,r}(s)-Z_{2,r}(s)},s\in[r,t]\}$ is a martingale. Note that $\langle X_1,X_2\rangle_s^c\equiv 0$. Then by It\^o's formula, we have
	\begin{eqnarray*}
		\ar\ar\!\!\!\!\!\!\!\! G(\mathbf u_{t,t},\mathbf X(t),\mathbf Z_r(t))=G(\mathbf u_{r,t},\mathbf X(r),\mathbf Z_r(r))+\sum_{i=1}^2\int_r^t G_i'(\mathbf u_{s-,t},\mathbf X(s-),\mathbf Z_r(s-))\d u_{i,s,t}\\
		\ar\ar+\sum_{i=1}^2\int_r^t G_{i+2}'(\mathbf u_{s-,t},\mathbf X(s-),\mathbf Z_r(s-))\d X_i(s)+\sum_{i=1}^2\int_r^t G_{i+4}'(\mathbf u_{s-,t},\mathbf X(s-),\mathbf Z_r(s-))\d Z_{i,r}(s)\\
		\ar\ar+\frac{1}{2}\sum_{i=1}^2\int_r^t G_{i+2,i+2}''(\mathbf u_{s-,t},\mathbf X(s-),\mathbf Z_r(s-))\d \langle X_i,X_i\rangle_s^c\\
		\ar\ar+\sum_{s\in(r,t]}[G(\mathbf u_{s,t},\mathbf X(s),\mathbf Z_r(s))-G(\mathbf u_{s-,t},\mathbf X(s-),\mathbf Z_r(s-))\\
		\ar\ar\quad\quad\quad\quad\quad -G_1'(\mathbf u_{s-,t},\mathbf X(s-),\mathbf Z_r(s-))\triangle u_{1,s,t}-G_2'(\mathbf u_{s-,t},\mathbf X(s-),\mathbf Z_r(s-))\triangle u_{2,s,t}\\
		\ar\ar\quad\quad\quad\quad\quad -G_3'(\mathbf u_{s-,t},\mathbf X(s-),\mathbf Z_r(s-))\triangle X_1(s)-G_4'(\mathbf u_{s-,t},\mathbf X(s-),\mathbf Z_r(s-))\triangle X_2(s)\\
		\ar\ar\quad\quad\quad\quad\quad -G_5'(\mathbf u_{s-,t},\mathbf X(s-),\mathbf Z_r(s-))\triangle Z_{1,r}(s)-G_6'(\mathbf u_{s-,t},\mathbf X(s-),\mathbf Z_r(s-))\triangle Z_{2,r}(s)]\\
		\ar=\ar G(\mathbf u_{r,t},\mathbf X(r),\mathbf Z_r(r))-\sum_{i=1}^2\int_r^t X_i(s-)G(\mathbf u_{s-,t},\mathbf X(s-),\mathbf Z_r(s-))\mathbf 1_{A^c}(s)\d u_{i,s,t}\\
		\ar\ar-\sum_{i=1}^2\int_r^t u_{i,s-,t}G(\mathbf u_{s-,t},\mathbf X(s-),\mathbf Z_r(s-))\mathbf 1_{A^c}(s)[\d X_s^{\rm mart.}-X_i(s-)b_{ii}(\d s)+X_j(s-)\bar b_{ji}(\d s)]\\
		\ar\ar-\sum_{i=1}^2\int_r^t G(\mathbf u_{s-,t},\mathbf X(s-),\mathbf Z_r(s-))\mathbf 1_{A^c}(s)X_i(s)\zeta_i(\d s)\\
		\ar\ar+\sum_{i=1}^2\int_r^t u_{i,s-,t}^2G(\mathbf u_{s-,t},\mathbf X(s-),\mathbf Z_r(s-))X_i(s-)c_i(\d s)\\
		\ar\ar+\sum_{s\in(r,t]\cap A^c}G(\mathbf u_{s-,t},\mathbf X(s-),\mathbf Z_r(s-))K(\mathbf u_{s-,t},\triangle\mathbf X(s))\\
		\ar\ar+\sum_{s\in(r,t]\cap A}[G(\mathbf u_{s,t},\mathbf X(s),\mathbf Z_r(s))-G(\mathbf u_{s-,t},\mathbf X(s-),\mathbf Z_r(s-))].
	\end{eqnarray*}
	By $\hat N_0(\d s,\d \mathbf z)$ in the proof of Theorem \ref{SIE's weak solution}, we have
	\begin{eqnarray*}
		&&\sum_{s\in(r,t]\cap A^c}G(\mathbf u_{s-,t},\mathbf X(s-),\mathbf Z_r(s-))K(\mathbf u_{s-,t},\triangle\mathbf X(s))=\textrm{mart.}\\
		&&\quad\quad\quad+\sum_{i=1}^2\int_r^t\int_{\mathbb R_+^2\backslash\{\mathbf 0\}} G(\mathbf u_{s-,t},\mathbf X(s-),\mathbf Z_r(s-)) K(\mathbf u_{s-,t},\mathbf z)\mathbf 1_{A^c}(s)X_i(s-)m_i(\d s,\d\mathbf z).
	\end{eqnarray*}
Since
 \begin{equation*}
\mathbf P_{r,\mathbf x}(G(\mathbf u_{t,t},\mathbf X(t),\mathbf Z_r(t)))=\mathbf P_{r,\mathbf x}(G(\mathbf u_{r,t},\mathbf X(r),\mathbf Z_r(r))),
 \end{equation*}
taking expectations on both sides of the above equation, we have, for $s\in A^c$,
 \begin{eqnarray*}
\d u_{i,s,t}=-\zeta_i(\d s)+u_{i,s,t}b_{ii}(\d s)-u_{j,s,t}\bar b_{ij}(\d s)+u_{i,s,t}^2 c_i(\d s)+\int_{\mathbb R_+^2\backslash\{\mathbf 0\}} K(\mathbf u_{s,t},\mathbf z)m_i(\d s,\d\mathbf z).
 \end{eqnarray*}
Combining above two equations, we have $\mathbf u_{r,t}$ satisfies \eqref{u_r,t}. \eproof

\medskip\noindent\textit{Proof of Theorem~\ref{positive integral functionals}.~} Let $\mathbf u_{r,t}$ be defined as in Proposition \ref{positive integral functionals prop} with $\bm\lambda=\mathbf 0$. Define $\mathbf w_{s,t}=\mathbf u_{s,t}+\bm\zeta(\{s\})$, where $\bm\zeta(\{s\})=(\zeta_1(\{s\}),\zeta_2(\{s\}))$. Then
 \beqnn
\mathbf P_{r,\mathbf x}\exp\bigg\{-\sum_{i=1}^2\int_{[r,t]} X_i(s)\zeta_i(\d s)\bigg\}\ar=\ar\mathbf P_{r,\mathbf x}\exp\bigg\{-\langle\mathbf X(r),\bm\zeta(\{r\})\rangle-\sum_{i=1}^2\int_r^t X_i(s)\zeta_i(\d s)\bigg\}\\
 \ar=\ar
\exp\bigg\{-\langle\mathbf x,\mathbf u_{r,t}+\bm\zeta(\{r\})\rangle\bigg\}=\exp\bigg\{-\langle\mathbf x,\mathbf w_{r,t}\rangle\bigg\},
 \eeqnn
where $r\mapsto\mathbf w_{r,t}$ is the unique positive solution to
 \beqnn
w_{i,r,t}\ar=\ar u_{i,r,t}+\zeta_i(\{r\})\\
 \ar=\ar
\zeta_i(\{r\})+\int_r^t\zeta_i(\d s)-\int_r^t \big[u_{i,s,t}+\zeta_i(\{s\})\big]b_{ii}(\d s)+\int_r^t \big[u_{j,s,t}+\zeta_j(\{s\})\big]b_{ij}(\d s)\\
 \ar\ar
-\int_r^t u_{i,s,t}^2 c_{i}(\d s)-\int_r^t\int_{\mathbb R_+^2\backslash\{\mathbf 0\}}K_i(\mathbf u_{s,t}+\bm\zeta(\{s\}),\mathbf z)m_i(\d s,\d \mathbf z)\\
 \ar=\ar
\int_{[r,t]}\zeta_i(\d s)-\int_r^t w_{i,s,t}b_{ii}(\d s)+\int_r^t w_{j,s,t}b_{ij}(\d s)-\int_r^t w_{i,s,t}^2 c_{i}(\d s)\\
 \ar\ar
-\int_r^t\int_{\mathbb R_+^2\backslash\{\mathbf 0\}}K_i(\mathbf w_{s,t},\mathbf z)m_i(\d s,\d \mathbf z).
 \eeqnn
That proves the result. \qed

\vskip 0.3cm
{\bf Acknowledgements} 
We are grateful to the Laboratory of Mathematics and Complex Systems (Ministry of Education) for providing us the research facilities. This research is supported by the National Key R\&D Program of China (No. 2020YFA0712901).

\end{document}